\documentclass[leqno]{amsart}
\usepackage{amsmath,amsfonts,amssymb}
\usepackage{amsthm}
\usepackage[dvips]{epsfig}

\def\e
{\varepsilon}

\def\C{{\mathbb C}}
\def\R{{\mathbb R}}
\def\N{{\mathbb N}}
\def\Z{{\mathbb Z}}

\def\O{{\mathcal O}}
\def\si{\sigma}
\def\u{{\tt u}}
\def\<{\langle}
\def\>{\rangle}
\def\DT{{\Delta t}}
\def\DX{{\Delta x}}
\def\({\left(}
\def\){\right)}

\def\Eq#1#2{\mathop{\sim}\limits_{#1\rightarrow#2}}
\def\Tend#1#2{\mathop{\longrightarrow}\limits_{#1\rightarrow#2}}

\def\d{{\partial}}
\numberwithin{equation}{section}

\theoremstyle{plain}
\newtheorem{theorem}{Theorem}[section]
\newtheorem{lemma}[theorem]{Lemma}

\newtheorem{proposition}[theorem]{Proposition}

\theoremstyle{definition}

\theoremstyle{remark}

\numberwithin{equation}{section}

\begin{document}
\title[Numerical aspects of NLS with caustics]{Numerical aspects of
  nonlinear Schr\"odinger equations in the presence of caustics}
\author[R. Carles]{R\'emi Carles}
\address[R. Carles]{Institut CNRS Pauli\\ Wolfgang
    Pauli Institute c/o Fak. f. Mathematik.\\
        Univ. Wien, UZA 4\\   
        Nordbergstr.~15\\ A-1090 Wien\\ Austria\footnote{Permanent
  address: Univ. Montpellier 2, UMR CNRS 5149, Math\'ematiques, CC
  051, Place Eug\`ene Bataillon, 34095 Montpellier cedex 5, France.}} 
\email{Remi.Carles@math.cnrs.fr}
\author[L. Gosse]{Laurent Gosse}
\address[L. Gosse]{Istituto
 per le Applicazioni del Calcolo (sezione di Bari)\\ Via G. Amendola
 122\\ 70126 Bari\\ Italy}
\email{l.gosse@ba.iac.cnr.it}
\thanks{This work was partially
  supported by o Centro de Matem\'atica e Aplica\c c\~oes Fundamentais
  (Lisbon), funded by FCT as contract POCTI-ISFL-1-209, and by the
  Austrian Ministry of Science via its  
grant for the Wolfgang Pauli Institute and by the Austrian Science
Foundation (FWF) via the START Project (Y-137-TEC).}
\begin{abstract}
The aim of this text is to develop on the asymptotics
  of some 1-D nonlinear Schr\"odinger equations from both the
  theoretical and the numerical perspectives, when a caustic is
  formed. We review rigorous 
  results in the field and give some heuristics in cases where
  justification is still needed. The scattering operator theory is
  recalled. Numerical 
  experiments are carried out on the focus point singularity for which
  several results have been proven rigorously. Furthermore, the scattering
  operator is numerically studied. Finally, experiments on the cusp
  caustic are displayed, and similarities with the focus point are
  discussed. 
\end{abstract}
\maketitle

\section{Introduction}
\vspace*{-0.5pt}
\noindent
We present a numerical study of the semi-classical solutions
to the following nonlinear Schr\"odinger equations with $\e \ll 1$,
\begin{equation}
  \label{eq:nls0}
  i\e \d_t \u^\e +\frac{\e^2}{2}\Delta \u^\e
  =|\u^\e|^{2\si}\u^\e,\quad(t,x)\in \R_+\times\R^n\quad  ;\quad
  \u^\e_{\mid t=0} = \e^p 
  f(x)e^{i\phi_0(x)/\e},
\end{equation}
when a caustic (a point or a cusp) is formed, that is to say, beyond
{\it breakup  
time}. Since the nonlinearity is homogeneous, the change of unknown function 
$u^\e = \e^{p}\u^\e$ shows that \eqref{eq:nls0} is equivalent to:
\begin{equation}
  \label{eq:nls1}
  i\e \d_t u^\e +\frac{\e^2}{2}\Delta u^\e
  =\e^{2\si p}|u^\e|^{2\si}u^\e\quad ;\quad u^\e_{\mid t=0} =
  f(x)e^{i\phi_0(x)/\e},
\end{equation}
so that we can always consider initial data of order $\O(1)$.

There are several motivations to study the behavior of \eqref{eq:nls1}
when a caustic is formed. First, on a purely academic level, we
recall that the description of the caustic crossing is complete in the
case of linear equations; see \cite{Du}. For nonlinear equations,
very interesting formal computations were proposed in \cite{HK87} (we
recall the main idea in Section~\ref{sec:analyt} below). For
\emph{dissipative} nonlinear wave equations, Joly, M\'etivier and Rauch
\cite{JMRTAMS95,JMRMemoir} have proved that the amplification of the
wave near the caustic can ignite the dissipation phenomenon in such a
way that the oscillations (that carry highest energy) are absorbed. 
The above nonlinear Schr\"odinger equation is the simplest model of a
\emph{conservative, nonlinear equation}. The mass and the energy of
the solution are independent of time (see \eqref{eq:conserv}
below). Therefore, different nonlinear mechanisms are expected. We
recall in Section~\ref{sec:analyt} some results that have been
established rigorously, and give heuristic arguments to extend these
results. This serves as a guideline for the numerical experiments
proposed after. 
\smallbreak

Second, \eqref{eq:nls1} may be considered as a
simplified model for Bose--Einstein condensation, which may be 
modeled (see e.g. \cite{DGPS,PiSt}) by:
\begin{equation}
  \label{eq:nlsharmo}
  i\e \d_t u^\e +\frac{\e^2}{2}\Delta u^\e
  =\omega^2 \frac{|x|^2}{2}u^\e +\e^{2}|u^\e|^{2\si}u^\e,
\end{equation}
with $\si =2$ if $n=1$, and $\si =1$ if $n=2$ or $3$. 
The power $\e^2$ in front of the nonlinearity depends on the r\'egime
considered, and in particular on the respective scales of different
parameters (see e.g. \cite{BurqZworski} and references therein). 
The role of the harmonic potential $|x|^2$ is to model a magnetic
trap. In the semi-classical limit $\e\to 0$ for the linear equation,
this potential causes focusing at the origin for solutions whose data are
independent of $\e$. This is to be compared with the case of
\eqref{eq:nls1} with initial quadratic oscillations as considered
below: the initial quadratic oscillations force the solution to
concentrate at one point in the limit $\e \to 0$. The parallel between
\eqref{eq:nls1} and \eqref{eq:nlsharmo} 
was extended and justified in \cite{CaIHP} for these nonlinear
equations.
\smallbreak

From both points of view, when a caustic point is formed, the caustic
crossing may be described in terms of the scattering operator
associated to 
\begin{equation*}
  i\d_t \psi +\frac{1}{2}\Delta\psi = |\psi|^{2\si}\psi.
\end{equation*}
This aspect is recalled in Section~\ref{sec:analyt}. For this reason,
we also pay a particular attention to this operator, independently of
the above semi-classical limit. Note  that besides the existence
of this operator, very few of its properties (dynamical, for instance)
are known. 
\smallbreak

In this paper, we
always assume $2\si p\ge 1$: one of the reasons is that when
$0\le 2\si p<1$, instability occurs, see
\cite{CaARMA,CaBKW}. Suppose for instance that $u^\e$ solves
\eqref{eq:nls1}, and that $\widetilde u^\e$ solves \eqref{eq:nls1},
where $f$ replaced by $(1+\delta^\e)f$, where $\delta^\e$ is a sequence
of real numbers going to zero as $\e \to 0$. Then there are some
choices of $\delta^\e$ for which 
\begin{equation*}
  \liminf_{\e \to 0}\|u^\e(t^\e)- \widetilde u^\e(t^\e)\|_{L^2}>0,
\end{equation*}
for some sequence of time $t^\e\to 0$ (see \cite{CaARMA}, and
\cite{BurqZworski} for a similar phenomenon with different initial
data). Therefore, 
producing reliable numerical tests in 
the case $0\le 2\si p<1$ (which is super-critical as far as WKB
analysis is concerned \cite{CaARMA,CaBKW}) seems to be a very delicate
issue, that we leave out in the present paper. 
\smallbreak

The rest of this paper is structured as follows. In
Section~\ref{sec:analyt}, we recall the general approach of WKB
analysis for the Schr\"odinger equation, the arguments of \cite{HK87},
and the rigorous results available for the semi-classical limit of
\eqref{eq:nls1} when a caustic reduced to a  point is formed. We then
recall the definition of the scattering operator. We also give
heuristic arguments to tackle the case of a ``supercritical focal
point'', and to guess what the critical indices are when a cusp
caustic is formed, instead of a focal point. 
In Section~\ref{sec:numgen}, we present the different strategies that
have been followed in the literature to study numerically the
semi-classical limit for nonlinear Schr\"odinger equations. 
Numerical experiments on the semi-classical limit for \eqref{eq:nls1}
in the presence of a focal point appear in Section~\ref{sec:foc}, and
the scattering operator is simulated in Section~\ref{sec:scattnum}. 
We present the numerical experiments of the semi-classical limit for
\eqref{eq:nls1} in the presence of a cusp caustic in
Section~\ref{sec:cusp}, and make conclusive remarks in
Section~\ref{sec:concl}. 

\section{Analytical approach}\label{sec:analyt}
\vspace*{-0.5pt}
\noindent
\subsection{Semi-classical limit of the free Schr\"odinger equation}
Consider the initial value problem, for $(t,x)\in \R_+\times\R^n$:
\begin{equation}
  \label{eq:schrodlibre}
  i\e \d_t v^\e +\frac{\e^2}{2}\Delta v^\e =0\quad ;\quad v^\e_{\mid
  t=0} = f(x)e^{i\phi_0(x)/\e}.
\end{equation}
The aim of WKB methods is to describe the asymptotic behavior of
$v^\e$ as $\e \to 0$. For instance, $\e$ can be related to the Planck
constant, and the asymptotic behavior of $v^\e$ is expected to yield a
good description of $v^\e$ when $\e$ is fixed, but small compared to
the other parameters. More precisely, seek
$v^\e$ of the form
\begin{equation}
  \label{eq:bkw}
  v^\e(t,x)\sim e^{i\phi(t,x)/\e}\(a_0(t,x)+\e a_1(t,x)+\ldots\)
  \quad\text{as }\e \to 0.
\end{equation}
Plugging this expansion into \eqref{eq:schrodlibre} and canceling the
$\O(\e^0)$ term, we see that
the phase $\phi$ must solve the eikonal equation:
\begin{equation}
  \label{eq:eikonale}
  \d_t \phi +\frac{1}{2}|\nabla \phi|^2 =0\quad ;\quad \phi_{\mid
  t=0}=\phi_0. 
\end{equation}
To cancel the $\O(\e^1)$ term,  the leading order amplitude
solves the transport equation:
\begin{equation}
  \label{eq:transport}
  \d_t a_0 +\nabla\phi \cdot \nabla a_0+\frac{1}{2}a_0 \Delta \phi
  =0\quad ;\quad a_{0\mid t=0}=f.
\end{equation}
The eikonal equation \eqref{eq:eikonale} is solved thanks to
Hamilton-Jacobi theory\footnote{but not according to the theory of viscosity
solutions! See e.g. \cite{koko}.}: $\phi$ is constructed locally in space and
time (see e.g. \cite{CaBKW} for a discussion on this aspect). Even if
$\phi_0$ is smooth, $\phi$ develops singularities in finite time in
general: the locus where $\phi$ is singular is called \emph{caustic}
(see e.g. the second volume of \cite{Hormander}). When $\phi$ becomes
singular, all the terms $a_0,a_1,\ldots$ may become singular as well. 
One easily observes that (\ref{eq:transport}) admits a ``divergence form":
$\d_t |a_0|^2+ \nabla \cdot (|a_0|^2 \nabla \phi)=0$.
To illustrate this general discussion, we consider two examples that
will organize the rest of this paper.\\

\noindent\emph{Example (Quadratic phase). }
  Let $\phi_0(x) = -\frac{|x|^2}{2}$. Then \eqref{eq:eikonale} and
  \eqref{eq:transport} can be  solved explicitly:
  \begin{equation*}
    \phi(t,x)=\frac{|x|^2}{2(t-1)}\quad ;\quad
    a(t,x)=\frac{1}{(1-t)^{n/2}}f\( \frac{x}{1-t}\).
  \end{equation*}
This shows that as $t\to 1$, $\phi$ and $a$ become singular: the wave
$u^\e$ focuses at the origin. This example can be viewed as the smooth
counterpart of the Cauchy problem
\begin{equation*}
  i\d_t \psi+\frac{1}{2}\Delta \psi=0\quad ;\quad \psi_{\mid t=0}=
  e^{-i\frac{|x|^2}{2}} .
\end{equation*}
Fourier analysis shows that $\psi_{\mid t=1}=\delta$, the Dirac
measure at the origin.\\

Of course, the solution of \eqref{eq:schrodlibre} can be represented
as an oscillatory integral:
\begin{equation}\label{eq:intosc}
  v^\e(t,x) = \frac{1}{(2i\pi t)^{n/2}}\int e^{i\frac{|x-y|^2}{2\e
  t} +i\frac{\phi_0(y)}{\e}}f(y)dy. 
\end{equation}
The caustic set is exactly the locus where the
critical points for the phase
\begin{equation*}
 \Phi_{t,x}(y)= \frac{|x-y|^2}{2 t}+\phi_0(y)
\end{equation*}
are degenerate. Outside the caustic, an approximation of $v^\e$ is
given by the stationary phase theorem (that we recalled as simply as 
possible in \cite{GosseWKB}). This leads us to the second
example we shall consider numerically:\\
\noindent
\emph{Example (Cusp). }
  Let $n=1$ and $\phi_0(x)=\cos x$. 
The set of degenerate critical points for $\Phi_{t,x}(y)$ (caustic) is
given implicitly by: 
\begin{equation*}
  {\mathcal C}=\left\{ (t,x)\in \R_+\times \R; \exists y\in \R ,\
  \frac{y-x}{t}=\sin y,\text{ and 
  }\frac{1}{t}=\cos y\right\}. 
\end{equation*}
As soon as $t\ge 1$, a caustic is formed (see Figure 2 in
\cite{GosseJinLi}). \\

When considering the asymptotic behavior of $u^\e$ beyond the caustic,
two main features must be considered: the creation of other phases\footnote{in
our mind, phases are always associated to oscillations whose
period depends on $\e$, and goes to infinity as $\e\to 0$ - rapid
oscillations. The wavelength may be proportional to $\e$, or, say, to
$\sqrt \e$.}, 
and  the Maslov index (see \cite{Du} for more general linear
equations). In the case of a focal point, the first aspect does not
exist: there is no creation of phase, and one phase is enough to
describe $v^\e$ past the focal point $(t,x)=(1,0)$. One can prove
easily the following result:
\begin{lemma}\label{lem1}
  Let $n \ge 1$ and $f\in {\mathcal S}(\R^n;\C)$. If
  $\phi_0(x)=-|x|^2/2$, then the asymptotic behavior (in $L^2(\R^n)$)
  of the solution $v^\e$ to \eqref{eq:schrodlibre} is given by:
  \begin{equation*}
    v^\e(t,x)\Eq \e 0
\left\{
  \begin{aligned}
    \frac{e^{i|x|^2/(2\e (t-1))}}{(1-t)^{n/2}}f\(\frac{x}{1-t}\)&
    \ \text{ if }t<1,\\
e^{-in\frac{\pi}{2}}\frac{e^{i|x|^2/(2\e
    (t-1))}}{(t-1)^{n/2}}f\(\frac{x}{1-t}\)& 
    \ \text{ if }t>1.
  \end{aligned}
\right.
  \end{equation*}
\end{lemma}
In this example, the Maslov index is $-n\pi/2$. In the case of the
cusp, three phases must be considered to describe the asymptotic
behavior of $v^\e$ beyond the caustic (see
e.g. \cite{GosseJinLi,GosseWKB}). \\
For future discussion on the numerical results, we state the following
more precise result, which follows from the stationary phase theorem:
\begin{lemma}\label{lem2}
  Let $n \ge 1$ and $f\in {\mathcal S}(\R^n;\C)$. If
  $\phi_0(x)=-|x|^2/2$, then the asymptotic behavior of the solution
  $v^\e$ to \eqref{eq:schrodlibre} at time $t=2$ is given by: 
  \begin{equation*}
    v^\e(2,x)=e^{-in\frac{\pi}{2}}e^{i|x|^2/(2\e)}f\(-x\)+ O(\e)\quad
    \text{in }L^2\cap  L^\infty(\R^n).
  \end{equation*}
\end{lemma}
\subsection{Caustics in the nonlinear case: heuristics}
\label{sec:heur}

Consider now the perturbation of \eqref{eq:schrodlibre} with a
nonlinear term:
\begin{equation}
  \label{eq:nls}
  i\e \d_t u^\e +\frac{\e^2}{2}\Delta u^\e
  =\e^{\alpha}|u^\e|^{2\si}u^\e\quad ;\quad u^\e_{\mid 
  t=0} = f(x)e^{i\phi_0(x)/\e}.
\end{equation}
The sign of the nonlinearity is chosen so that no finite time blow-up
occurs. The following two important
quantities are formally independent of time:
\begin{equation}
  \label{eq:conserv}
  \begin{aligned}
  \text{Mass: }& \|u^\e(t)\|_{L^2}=\text{Const.}= \|f\|_{L^2}.\\
\text{Energy: }& E^\e(t):= \frac{1}{2}\|\e \nabla u^\e(t)\|^2_{L^2}
+\frac{\e^{\alpha}}{\si +
  1}\|u^\e(t)\|_{L^{2\si+2}}^{2\si+2}=E^\e(0).
\end{aligned}
\end{equation}
We refer to \cite{CazCourant} for a justification. Fix the power
$2\si>0$ of the nonlinearity, and consider different values for
$\alpha$. Two notions of criticality arise: for the WKB methods on the
one hand, and for the caustic crossing on the other hand. This
discussion is presented in \cite{HK87} for conservation laws, and we
summarize it  in the case of \eqref{eq:nls}. Plugging an expansion of
the form \eqref{eq:bkw} into \eqref{eq:nls}, we see that the value
$\alpha =1$ is critical for the WKB methods: if $\alpha >1$, then the
nonlinearity does 
not affect the transport equation \eqref{eq:transport} (``linear
propagation''), while if 
$\alpha =1$, then the nonlinearity appears in the right hand side of
\eqref{eq:transport} (``nonlinear propagation''). Recall that in this
paper, we always assume $\alpha \ge 1$. Therefore, the eikonal equation
\eqref{eq:eikonale} is not altered: 
the geometry of the propagation remains the same as in the linear 
WKB approach,
and we have to face the same caustic sets. The idea presented in
\cite{HK87} consists in saying that according to the geometry of the
caustic, different notions of criticality exist, as far as $\alpha$ is
concerned, near the
caustic. In the linear setting \eqref{eq:schrodlibre}, the influence
of the caustic is relevant only in a neighborhood of this set
(essentially, in a boundary layer whose size depends on $\e$ and the
geometry of $\mathcal C$). View the nonlinearity in \eqref{eq:nls}
as a potential, and assume that the nonlinear effects are negligible
near the caustic: then $u^\e\sim v^\e$ near $\mathcal C$. 
View the term $\e^\alpha |u^\e|^{2\si}$ as a (nonlinear)
potential. 
The average nonlinear effect near $\mathcal C$ is expected to be:
\begin{equation*}
 \e^{-1}\int_{{\mathcal C}^\e}\e^\alpha |u^\e|^{2\si}\sim
 \e^{-1}\int_{{\mathcal C}^\e}\e^\alpha |v^\e|^{2\si},  
\end{equation*}
where ${\mathcal C}^\e$ is the region where caustic effects are
relevant, and the factor $\e^{-1}$ is due to the integration in
time (recall that there is an $\e$ in front of the time derivative in
\eqref{eq:nls}).  The idea of this heuristic argument is that when the
nonlinear 
effects are negligible near $\mathcal C$ (in the sense that the
uniform norm of $u^\e-v^\e$ is small compared to that of $v^\e$ near
${\mathcal C}^\e$), the above approximation should be valid. On the
other hand, it is expected that it ceases to be valid precisely when
nonlinear effects can no longer be neglected near the caustic:
$u^\e-v^\e$ is of the same order of magnitude as $v^\e$ in
$L^\infty({\mathcal C}^\e)$, or even larger. 

Practically, assume that in the linear case, $v^\e$ has an
amplitude $\e^{-\ell}$ in a boundary layer of size $\e^{k}$; then the
above quantity is
\begin{equation*}
  \e^{-1}\int_{{\mathcal C}^\e}\e^\alpha |v^\e|^{2\si} \sim
  \e^{-1}\e^\alpha | \e^{-\ell}|^{2\si}\e^{k}.
\end{equation*}
The value $\alpha $ is then critical when the above cumulated effects
are not negligible:
\begin{equation*}
  \alpha_c = 1+2\ell \si -k.
\end{equation*}
When $\alpha>\alpha_c$, the nonlinear effects are expected to be
negligible near the caustic: resuming the terminology of \cite{HK87},
we speak of ``linear caustic''. The case $\alpha = \alpha_c$ is called
``nonlinear caustic''.
To conclude this paragraph, we examine this approach in the case of
our two examples.
In the case of a focal point, we have $k=1$ and $\ell=n/2$. This leads us
to the value:
\begin{equation*}
  \alpha_c(\text{focal point}) = n\si.
\end{equation*}
In the case of the cusp in dimension one, we have $k=2/3$ and
$\ell=1/3$ (which can be viewed thanks to the Airy function and its
asymptotic expansion, see e.g. \cite{Du,Hormander,HK87} or
\cite{Ludwig}), which yields:  
\begin{equation*}
  \alpha_c(\text{cusp in 1D}) = \frac{2\si +1}{3}.
\end{equation*}
One aspect of the numerical experiments presented below is to test
this notion of criticality in those two examples.

\subsection{Justification for a focal point, and more heuristics}
\label{sec:justif}
In this paragraph, we assume $\phi_0(x)= -|x|^2/2$. A complete
justification of the above discussion is available \cite{CaIUMJ}:\\
\smallbreak
\begin{center}
\begin{tabular}[c]{l|c|c}
 &$\alpha >n\si$  & $\alpha =n\si$  \\
\hline
$\alpha >1$ &Linear caustic, & Nonlinear caustic,\\
&linear propagation   & linear propagation \\
\hline
$\alpha =1$&Linear caustic, &Nonlinear caustic,\\
&nonlinear propagation   & nonlinear propagation \\
\end{tabular}
\end{center}
\smallbreak

Consider $t\in [0,2]$, which includes the caustic crossing. The above
tables means: 
\begin{itemize}
\item If $\alpha >\max (1,n\si)$, then $u^\e$ can be approximated
  by $v^\e$ for $t\in [0,2]$. 
\item If $\alpha =1>n\si$, then the nonlinearity is negligible
  near the focal point, but not away from it. 
\item If $\alpha =n\si>1$, then nonlinear effects are relevant
  near the focal point, and only near the focal point.
\item If $\alpha =n\si=1$, then the nonlinearity is never
  negligible. 
\end{itemize}
We give some precisions in some cases of interest for the
numerics presented below. \\

In the one-dimensional case $n=1$, the following pointwise estimate is
proved in \cite{CaIUMJ} when $\alpha >\max (1,\si)$ or $\alpha =\si>1$:
\begin{equation}\label{eq:estponctuelle}
  |u^\e(t,x)|\le \frac{C}{\sqrt{|t-1|+\e}}\cdot
\end{equation}
Setting $w^\e=u^\e-v^\e$, we see that
\begin{equation*}
  i\e \d_t w^\e +\frac{\e^2}{2}\d_x^2 w^\e = \e^\alpha
  |u^\e|^{2\si}u^\e \quad ;\quad w^\e_{\mid t=0}=0.
\end{equation*}
The usual energy estimate yields, for $t\ge 0$:
\begin{equation*}
  \|w^\e(t)\|_{L^2}\le \e^{-1}\int_0^t \e^\alpha
  \big\| |u^\e(\tau)|^{2\si}u^\e(\tau)\big\|_{L^2}d\tau.
\end{equation*}
Using \eqref{eq:estponctuelle}, we infer, for $t\in [0,2]$:
\begin{equation*}
  \|w^\e(t)\|_{L^2}\lesssim
  \e^{\alpha-1}\(\int_{\{|\tau-1|>\e\}\cap \{\tau \in
  [0,2]\}}\frac{d\tau}{|\tau 
  -1|^\si}+\int_{|\tau-1|\le \e}\frac{d\tau}{\e^\si} \)\lesssim
  \e^{\alpha -\si}. 
\end{equation*}
Using the operator $\e\d_x$ and $x/\e+i(t-1)\d_x$, and
Gagliardo--Nirenberg inequalities as in \cite{CaIUMJ}, we find:
\begin{lemma}\label{lem3}
  Let $n=1$, $\alpha >\max (1,\si)$, $f\in {\mathcal S}(\R^n;\C)$, and
  $\phi_0(x)=-|x|^2/2$. Then we have, for the solutions of
  \eqref{eq:schrodlibre} and \eqref{eq:nls}:
  \begin{equation*}
    \sup_{0\le t\le 2} \big\| u^\e(t)-v^\e(t)\big\|_{L^2}\le C
    \e^{\alpha -\si}\quad ;\quad \big\| u^\e(t)-v^\e(t)\big\|_{L^\infty}\le C
    \frac{\e^{\alpha -\si} }{\sqrt{|t-1|+\e}}, \quad t\in [0,2].
  \end{equation*}
In particular, Lemma~\ref{lem3} implies, at time $t=2$:
\begin{equation*}
    u^\e(2,x)=e^{-in\frac{\pi}{2}}e^{i|x|^2/(2\e)}f\(-x\)+
    O\(\e^{\min(1,\alpha -\si)}\)\quad
    \text{in }L^2\cap  L^\infty(\R).
  \end{equation*}
\end{lemma}
(The above result is true also when $\alpha=\si>1$, but becomes far
less interesting.)
We now explain the critical case $\alpha =n\si>1$. The nonlinear effects near
the focal point are described in terms of the scattering operator
associated to the nonlinear Schr\"odinger equation. We rapidly present
this operator $S$ in Section~\ref{sec:scatt}. We have then:
\begin{equation*}
    u^\e(t,x)\Eq \e 0 
\left\{
  \begin{aligned}
    \frac{e^{i|x|^2/(2\e (t-1))}}{(1-t)^{n/2}}f\(\frac{x}{1-t}\)&
    \ \text{ if }t<1,\\
e^{-in\frac{\pi}{2}}\frac{e^{i|x|^2/(2\e
    (t-1))}}{(t-1)^{n/2}}Zf\(\frac{x}{1-t}\)& 
    \ \text{ if }t>1,
  \end{aligned}
\right.
  \end{equation*}
where $Z= {\mathcal F}\circ S\circ {\mathcal F}^{-1}$ is the conjugate
of $S$ by the Fourier transform (see 
\cite{CaIUMJ}),
\begin{equation*}
  {\mathcal F} \varphi(\xi)=\frac{1}{(2i\pi)^{n/2}}\int e^{-ix\cdot
  \xi}\varphi(x)dx. 
\end{equation*}
(Since $S$ is a nonlinear operator, the normalization of the Fourier
transform is an important detail.)
\smallbreak 

To conclude this paragraph, we give a few hints of what happens or may
happen when the propagation is linear, and the caustic is
super-critical, that is $1<\alpha <n\si$. First, the conservation
of mass and energy seem to rule out the possibility of a concentration
of the form
\begin{equation}\label{eq:concentration}
  u^\e(1,x)\sim \frac{1}{\e^{n/2}}\varphi \(\frac{x}{\e}\),  
\end{equation}
for some function $\varphi$ independent of $\e$. The above
relation holds for $v^\e$ and $u^\e$ when $\alpha >\max(1,n\si)$,
with $\varphi ={\mathcal F}f$ (for $v^\e$, this is obvious from
\eqref{eq:intosc}). When $\alpha =n\si>1$ (linear propagation and
nonlinear focal point), the above relation still holds, with a
different profile $\varphi$ (see \cite{CaIUMJ}). Now we see that the
energy is bounded as $\e\to 0$:
\begin{equation*}
  E^\e(t)= E^\e(0)\Eq \e 0 \frac{1}{2}\|x f\|_{L^2}^2.
\end{equation*}
Plugging a concentrating profile as in \eqref{eq:concentration} in the
second term of the energy would yield, thanks to the conservation of mass:
\begin{equation*}
  \e^\alpha \|u^\e(1)\|_{L^{2\si+2}}^{2\si+2}\thickapprox \e^{\alpha
  -n\si}\Tend \e 0 +\infty,
\end{equation*}
which is impossible since the energy is the sum of two positive
terms. This suggests two possible effects: near $t=1$,  the
amplification of the solution (in terms of powers of $\e$) may be
weaker than in 
the linear case; on the other hand,  nonlinear effects near the 
caustic should affect the \emph{phase} of the solution $u^\e$ in a
rather strong way, causing the appearance of new frequencies. A partial
justification of the last assertion may be found in
\cite{CaCascade}. 
\smallbreak

Extrapolating this argument, we expect that in the supercritical case
for a cusp (with $n=1$: $\frac{2\si+1}{3}>\alpha>1$), new
frequencies appear. If as in 
the linear case,
three phases are necessary to describe the solution
past the caustic, then the nonlinear interaction of these phases might
reveal the presence of new frequencies, even on the modulus of $u^\e$
(see Sect.~\ref{sec:cusp} for numerical 
tests that seem to confirm this heuristics). 

\subsection{The scattering operator for NLS}
\label{sec:scatt}
To explain what the operator $S$ mentioned in the previous section is,
consider the nonlinear Schr\"odinger equation 
\begin{align}
  i\d_t \psi +\frac{1}{2}\Delta \psi & = |\psi|^{2\si}\psi\quad ;\quad
  (t,x)\in \R\times \R^n,\label{eq:NLS} \\
  U(-t)\psi(t)_{\mid t=t_0}&= \psi_-,\label{ci:NLS}
\end{align}
where $U(t)= e^{i\frac{t}{2}\Delta}$ is the propagator of the linear
equation. To construct the scattering operator, we first want to give a
meaning to \eqref{eq:NLS}--\eqref{ci:NLS} when $t_0=-\infty$. This
means that the nonlinear effects are asymptotically negligible as $t\to
-\infty$: for instance, we expect at least
\begin{equation*}
  \| U(-t)\psi(t) -\psi_-\|_{L^2}= \| \psi(t) -U(t)\psi_-\|_{L^2}\Tend
  t {-\infty} 0.
\end{equation*}
This gives a rigorous meaning to the relation $\psi(t,x)\sim
U(t)\psi_-(x)$ which aims at saying that  as time goes to $-\infty$,
the nonlinear dynamics associated to \eqref{eq:NLS} can be compared to
the free dynamics given by $e^{i\frac{t}{2}\Delta}$. 

To define a scattering operator, we want to be able to say that as
$t\to +\infty$ as well, the nonlinear effect are asymptotically
negligible. That is, there exists $\psi_+\in L^2(\R^n)$ such that
\begin{equation*}
  \| U(-t)\psi(t) -\psi_+\|_{L^2}= \| \psi(t) -U(t)\psi_+\|_{L^2}\Tend
  t {+\infty} 0.
\end{equation*}
The scattering operator is then defined as $S:\psi_-\mapsto
\psi_+$. Since our numerical experiments concern the one-dimensional
case, we recall the existence of the scattering operator in this
setting, and refer to
\cite{CazCourant,CW92,GOV94,GV79Scatt,NakanishiOzawa} for some
extensions to the multidimensional framework: define 
\begin{equation*}
  \Sigma := H^1\cap {\mathcal F}(H^1)=\{ f\in L^2(\R^n)\ ;\
  \|f\|_\Sigma := \|f\|_{L^2}+\|x f\|_{L^2} + \|\nabla
  f\|_{L^2}<\infty\}. 
\end{equation*}
\begin{proposition}[Scattering theory]\label{prop:scattering}
Let $n=1$, and assume $\si \geq \frac{1+\sqrt{17}}{4} (>1)$.
\begin{itemize}
\item For every  $\psi_- \in \Sigma $,
 there exists a unique $\varphi\in \Sigma$ 
 such that the maximal solution $\psi\in C(\R,\Sigma)$ to
 \eqref{eq:NLS} satisfies $\psi_{\mid t=0}=\varphi$ and 
$$\left\| U(-t)\psi(t)-\psi_-\right\|_\Sigma \Tend
t {-\infty} 0.$$
\item For every $\varphi\in \Sigma$, there exists a unique
 $\psi_+\in \Sigma$  
 such that the maximal solution $\psi\in C(\R,\Sigma)$ to
 \eqref{eq:NLS} with $\psi_{\mid t=0}=\varphi$ satisfies
$$\left\| U(-t)\psi(t)-\psi_+\right\|_\Sigma \Tend
t {+\infty} 0.$$
\end{itemize}
The scattering operator is $S:\psi_-\mapsto \psi_+$. When $\si >1$,
the above conclusions remain, in a neighborhood of the origin. When
$\si>1$, we also have:
\begin{itemize}
\item For every $\psi_-\in\Sigma$, there exist a unique solution
 $\psi\in C(\R,H^1)$ to 
 \eqref{eq:NLS} and a unique $\psi_+\in
  H^1(\R)$  such that:
 \begin{equation*}
   \left\| U(-t)\psi(t)-\psi_-\right\|_{\Sigma} \Tend
t {-\infty} 0\quad ;\quad \left\| U(-t)\psi(t)-\psi_+\right\|_{L^2} \Tend
t {+\infty} 0.
 \end{equation*}
\end{itemize}
\end{proposition}

When $\si \le 1$, the above conclusions are false: if $\si =1$ for
instance, and if $\psi_-\in L^2$ with
$U_0(-t)\psi(t)-\psi_- \to 0$ in $L^2$ as $t\to -\infty$, then
$\psi=\psi_-=0$. One cannot compare the nonlinear dynamics with the
free dynamics (see \cite{Barab,Strauss74,Strauss81,Ginibre}). \\
Note that even though the scattering is proven to exist, very few of its
features are known. We refer to \cite{CazCourant} for some algebraic
properties. At least, this operator is not trivial: near
the origin, it is a non-trivial perturbation of the identity (see
\cite{CaWigner}). 

\section{Numerical approximation of semi-classical Schr\"odinger
  equations}\label{sec:numgen} 
\vspace*{-0.5pt}
\noindent

Hereafter we restrict our discussion to the one-dimensional case, that
is to say $n=1$ in all the preceding considerations. 

\subsection{Rigorous results for general time-splitting schemes}

It is interesting to notice that in the numerical literature, the
nonlinear equation \eqref{eq:nls} is treated exactly the same way the
linear one \eqref{eq:schrodlibre} would be in the presence of a
potential $V$ on its right-hand side. The strategy is called {\it
  time-splitting}, in its first or second order version (Lie or Strang
splitting, see e.g. \cite{BBD}) where one alternates every time step
$\DT>0$ between the solving of the Laplace operator and the handling
of the (nonlinear) differential equation. According to
Section \ref{sec:scatt}, $U$ will still stand for the free propagator,
whereas we shall use $V$ as the ODE solver; Lie time-splitting
algorithms generate the following type of approximation for
\eqref{eq:nls}, 
$$u^\e(n\DT,.) \simeq u^\e_\DT(n\DT,.):=[V(\DT) \circ U(\DT)]^n u^\e(t=0,.),$$
and Strang splittings, 
$$u^\e(n\DT,.) \simeq \tilde u^\e_\DT(n\DT,.):=V(\DT/2) \circ U(\DT)
[V(\DT) \circ U(\DT)]^{n-1} V(\DT/2) u^\e(t=0,.).$$ 
Many references exist; let us quote only \cite{BBD,MPP,BJM1,BJM2}.

On the contrary, few rigorous convergence results are available, hence
we shall mainly recall the results from \cite{BBD} which quantify
accurately the {\it splitting errors} assuming each time-step is
performed exactly\footnote{but we shall see in the sequel that this is
  far from being the case!}. Under this assumption, there holds: 
\begin{proposition}(\cite{BBD}, Theorem 4.1) 
  For any $T>0$ and $u^\e(t=0,.) \in H^2$, there
  exists a constant $C$ depending on the initial data for
  \eqref{eq:nls} and $h_0$, such that for $\DT \in [0,h_0]$ and $n\DT
  <T$,  
$$\|u^\e(n\DT,.) -u^\e_\DT(n\DT,.)\|_{L^2}\leq C h_0.$$
If moreover $u^\e(t=0,.) \in H^4$, then there holds under the same assumptions:
$$\|u^\e(n\DT,.) -\tilde u^\e_\DT(n\DT,.)\|_{L^2}\leq C h^2_0.$$
\end{proposition}

This result is concerned with splitting errors only and relies on the
knowledge of the exact solution operators $U$ and $V$. In order to
stick to this framework in the context of smooth solutions, it is
rather natural to approximate $U$ by means of a Fourier scheme taking
advantage of optimized FFT routines, as proposed in the paper
\cite{BJM2}. Moreover, this will guarantee that the $L^2$ norm (Mass)
of the numerical solution will be conserved up to round-off
errors. Unfortunately, the Hamiltonian $E^\e(t)$ is generally not
preserved; a method conserving both quantities exists (see the
so--called MCN algorithm, page 253 of \cite{DSS}) but it wouldn't be
efficient in the semiclassical regime because of the results in
\cite{MPP}. 

\subsection{Specific issues with finite-difference discretizations}

This is the main purpose of the paper \cite{MPP} to illustrate the
(surprising) fact that in semi-classical regime, usual
finite-difference schemes for \eqref{eq:schrodlibre} can deliver very
wrong approximations without any particular sign of instability in
case very restrictive meshing constraints turn out to be
bypassed. This can be quite easily checked through the location of
caustics, for instance. The analysis of those standard schemes has
been carried out by means of Wigner measures, so the conclusions hold
essentially for the quadratic observables coming out of the wave
function itself. 

\subsection{The case of FFT-based schemes}

This class of schemes became popular after the publication of
\cite{BJM1,BJM2}, mainly because treating the differential part of
(\ref{eq:nls}) by means of a discrete Fourier transform looked very
much like being the best possible compromise in terms of meshing
constraints. Indeed, in the linear case where (\ref{eq:schrodlibre})
is supplemented with a potential $V(x)$ on its right-hand side, it was
shown that the time-step $\DT$ could be chosen independent of $\e$
whereas the space discretization has to satisfy $\DX=\O(\e)$. This was
already much better when compared to finite-differences; moreover, the
method is naturally $L^2$-conservative. In \cite{BJM2}, these authors
extended their ``Fourier framework" to the weakly nonlinear
Schr\"odinger equations of the form (\ref{eq:nls}).  
%
%
However, and despite the fact we do believe these 
``FFT time-split schemes'' realize the best numerical strategy in
terms of gridding, we 
shall point out some shortcomings of the method in the next section. 

\subsection{A rigorous framework for FFT-based schemes}
\label{sec:rig_error}

We present here a preliminary result about truncation errors 
in Lebesgue spaces for Fourier schemes; its proof follows directly
from the Strichartz estimates on the torus due to
J.~Bourgain \cite{bou} (see also \cite{bou_lect}),  and from the study
of FFT by M.~Taylor 
\cite[pp.~250--254]{taylor_PDE}. Its derivation is 
not obvious though as it applies directly to widely-used schemes like
the one recalled in the forthcoming section. We restrict our attention to  
the 1D free Schr\"odinger equation \eqref{eq:schrodlibre} with $\e=1$,
and periodic boundary  
conditions: $x \in {\mathbb T}:={\mathbb R}/2\pi{\mathbb Z}$.

Hence we start from 
$$
i\d_t \psi + \frac 1 2 \d_x^2 \psi=0, \qquad \psi(t=0,\cdot)=\zeta=
\sum_{j \in \Z} \hat \zeta_j e^{ijx},\quad 
x \in [0,2\pi]. 
$$
We have explicitly:
$$
\psi(t,x)=\sum_{j \in \Z} \hat \zeta_j e^{ij(x-jt/2)}.
$$
In order to investigate the behavior of the FFT-scheme involving a
finite even number $N \in 2\N$ of modes, we introduce the Discrete
Fourier Transform of a continuous function $f$ on $[0,2\pi]$ as follows:
$$
f^\#_k \stackrel{\text{def}}{=}\frac 1 N \sum_{j=1}^N f(2j\pi
/N)e^{-i 2\pi jk/N}. 
$$
From \cite[p.~252]{taylor_PDE}, we recall:
\begin{lemma}\label{lem_tay}If the continuous function $f$ has
a convergent Fourier series, then:
$$\forall k, \qquad f^\#_k = \sum_{j \in \Z} \hat f_{k +jN}.$$
\end{lemma}
Now, what the Fourier numerical scheme really computes is:
$$
\forall k \in \left\{-\frac N 2, \frac N 2\right \},
\qquad \hat \psi^\#_k(t)\stackrel{\text{def}}{=}e^{-i k^2 t/2}\hat
\zeta^\#_k. 
$$
Thus the corresponding numerical solution $\psi^{Num}$ is built:
$$
\psi^{Num}(t,x)=\sum_{k=-\frac N 2}^{\frac N 2}e^{i k(x-
kt/2)}\sum_{j \in \Z}\hat \zeta_{k +jN}. 
$$
The question is therefore to study the discrepancy between $\psi$ and
$\psi^{Num}$: 
\begin{align*}
\psi^{Num}(t,x)-\psi(t,x)&= \sum_{k=-\frac N 2}^{\frac N 2}e^{i
k(x- kt/2)}\sum_{j \in \Z}\hat \zeta_{k +jN} 
-\sum_{k \in \Z} \hat \zeta_k e^{ik(x-kt/2)}\\
&= \sum_{j \in \Z}\left(\sum_{k=-\frac N 2}^{\frac N 2}e^{i k(x-
  kt/2)}\hat \zeta_{k +jN} \right)- 
\sum_{k \in \Z} \hat \zeta_k e^{ik(x-kt/2)}\\
&=\sum_{j \in \Z^\star}\left(\sum_{k=-\frac N 2}^{\frac N 2}e^{i
  k(x- kt/2)}\hat \zeta_{k +jN} \right)+ 
\sum_{|k| > \frac N 2} \hat \zeta_k e^{ik(x-kt/2)}\\
&\stackrel{\text{def}}{=} I + II.
\end{align*}
The second term $II$ in this last equality can be interpreted as a
usual truncation error; it means that 
some high frequencies in the solution are lost when discretizing the
problem on a finite grid. From \cite[p.~253]{taylor_PDE}, we know that
the modulus of this term can be controlled as follows: 
$$
\sum_{|k| > \frac N 2} |\hat \zeta_k| \leq \frac 1 {N^m} \| \zeta \|_{C^{1+m}},
$$
which is satisfactory provided the initial data $\zeta$ is smooth,
e.g. $C^2(\R)$. However, the first term $I$ is far more delicate and
reveals the propagation of the errors coming from the DFT/FFT
itself. It can be controlled thanks to the  periodic Strichartz
estimates proved by J.~Bourgain. From \cite[pp.~16--17]{bou_lect}, we
recall: 
\begin{lemma}\label{lem_bou}
Given any complex-valued sequence $(a_k)_k$, 
the following estimates hold for the corresponding 
functions in $[0,2\pi] \times [0,2\pi]$:
\begin{align*}
&\Big\|\sum_k a_k e^{i(kx-k^2t)}\Big\|^2_{L^4({\mathbb T}\times
  {\mathbb T}) } \le C \sum_k |a_k|^2,\\  
&\Big\|\sum_{|k| \leq N} a_k e^{i(kx-k^2t)}\Big\|^2_{L^6({\mathbb
    T}\times {\mathbb T})} \le C e^{2\log 
N/\log \log N}\sum_k |a_k|^2. 
\end{align*}
\end{lemma}
The second estimate looks more attractive as it ``sees" the finite
number of modes. We therefore  deduce that the first term can be
controlled in $L^6$ by means of: 
$$
C \exp\left(\frac{\log N/2}{\log \log N/2}\right)
\left (\sum_{j \in \Z^\star, |k| \leq \frac N 2} |\hat
  \zeta_{k+jN}|^2\right)^{\frac 1 2}. 
$$
For instance, if $\zeta$ is a finite superposition of Fourier modes,
then it is clear that this 
term cancels for $N$ large enough as $j \not = 0$ in the summation;
obviously, the second term $II$ vanishes too. 


The general case isn't completely clear yet.

\section{Experiments on the focal point}\label{sec:foc}
\vspace*{-0.5pt}
\noindent

This section aims at visualizing the asymptotics previously recalled;  
namely we shall compare numerical approximations of \eqref{eq:nls} and
\eqref{eq:schrodlibre} 
in 1D ($n=1$) for various values of the parameters $\alpha$ and
$\sigma$. The initial wave function 
is rather simple:
$$
v^\e(t=0,x)=u^\e(t=0,x)=\exp \left(-\Big(2+\frac i {2\e}\Big)(x-\pi)^2\right), \qquad x\in [0,2\pi].
$$
Numerical results have been obtained through the time-splitting FFT
schemes recalled in the 
previous section; we used 1024 modes and fixed $\e=1/150$. It is
convenient to observe results 
in $t=2$ since $|v^\e(t=2,.)|=|v^\e(t=0,.)|$.

\subsection{Subcritical case}
\label{sec:foc_free}

\begin{figure}[ht]
\centerline{\epsfig{file=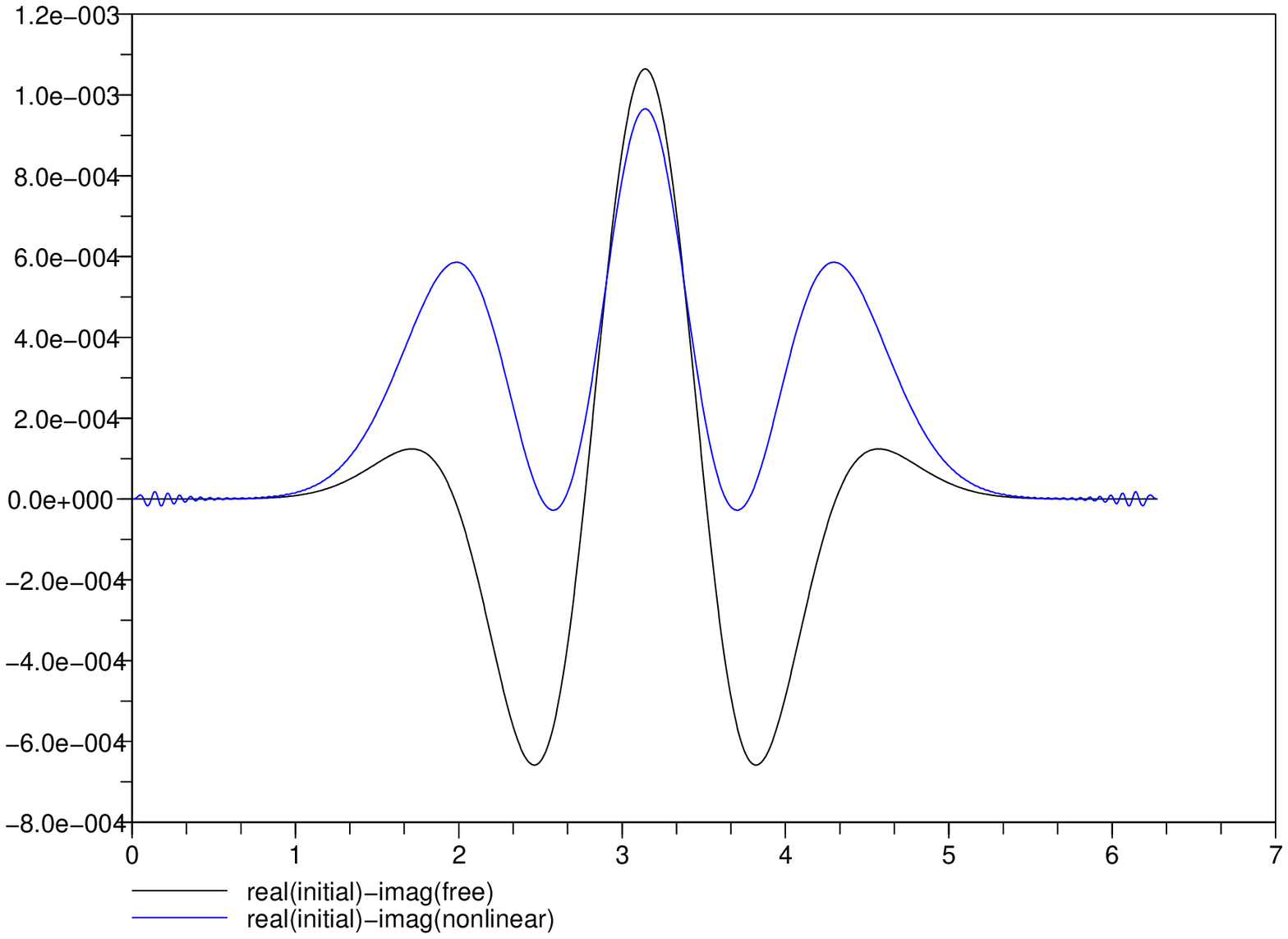,width=0.45\linewidth}
\epsfig{file=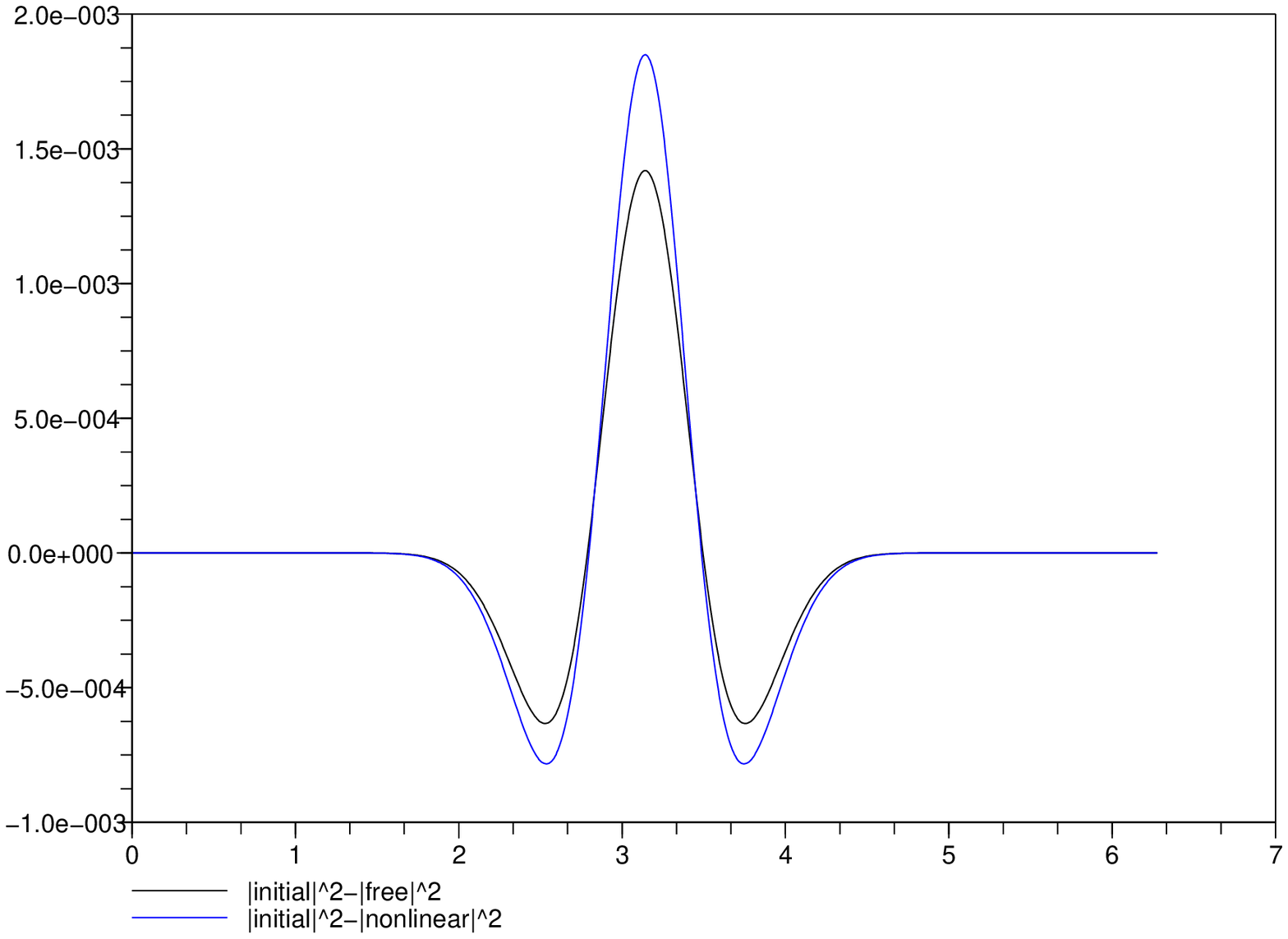,width=0.45\linewidth}}
\caption{Absolute errors on the wave functions (left) and on the
  modulus (right) at  
$T=2$ with $\sigma=2$ and $\alpha=2.5$. (Subcritical case)} \label{pf1}
\end{figure}

This case corresponds to $\sigma=2$ and $\alpha=2.5$; we expect to
observe a decay of the absolute 
errors between $v^\e$ and $u^\e$ for values $\e \ll 1$. This is indeed
the case, but Fig.~\ref{pf1} 
shows even a bit more, namely it compares pointwise the following
quantities: (recall Lemma \ref{lem1}) 
$$\left\{\begin{array}{l}
\Re\left(v^\e(t=0,x)\exp(-i(x-\pi)^2/2\e)\right), \\
\Im\left(v^\e(t=2,x)\exp(i(x-\pi)^2/2\e)\right), \\
\Im\left(u^\e(t=2,x)\exp(i(x-\pi)^2/2\e)\right).
\end{array}\right.
$$
On the left in Fig.~\ref{pf1}, we obviously observe that the absolute
errors are slightly bigger when 
considering the solution of the nonlinear equation \eqref{eq:nls}, $u^\e$. 
However, even for the free solution $v^\e$, one sees that the error
doesn't vanish despite the fact 
no time-splitting algorithm is needed. As the way of discretizing the solution
reveals itself important, we include
here the corresponding {\sc Scilab} routine for the free equation:

{\tt clear;deff('[y]=phase(x)',['y=-0.5*(x-$\pi$)$.^2$';])

deff('[y]=position(x)',['y=exp(-2*(x-$\pi$)$.^2$)'])

deff('[y]=Az(x)',['y=position(x).*exp(i*phase(x)./epsilon)'])

NMAX=2$^{10}$;n=-(NMAX)/2:(NMAX/2)-1;epsilon=1.0/150;

XSTART=0;XSTOP=2*$\pi$;DX=(XSTOP-XSTART)/NMAX;XSTOP=XSTOP-DX;

a=XSTART:DX:XSTOP,initialdata=Az(a);

vepsilon=fftshift(fft(initialdata,-1));

vepsilon=exp(-i*epsilon*($n.^2$)).*vepsilon;

vepsilon=fft(fftshift(vepsilon),1);
%
%
}\\
Clearly, its outcome is in agreement with Lemma
\ref{lem1} since $\e$ is already quite small. The Maslov index is visible, 
up to an error around $10^{-3}$ for 1024 Fourier modes.

\subsection{Critical case}

\begin{figure}[ht]
\centerline{\epsfig{file=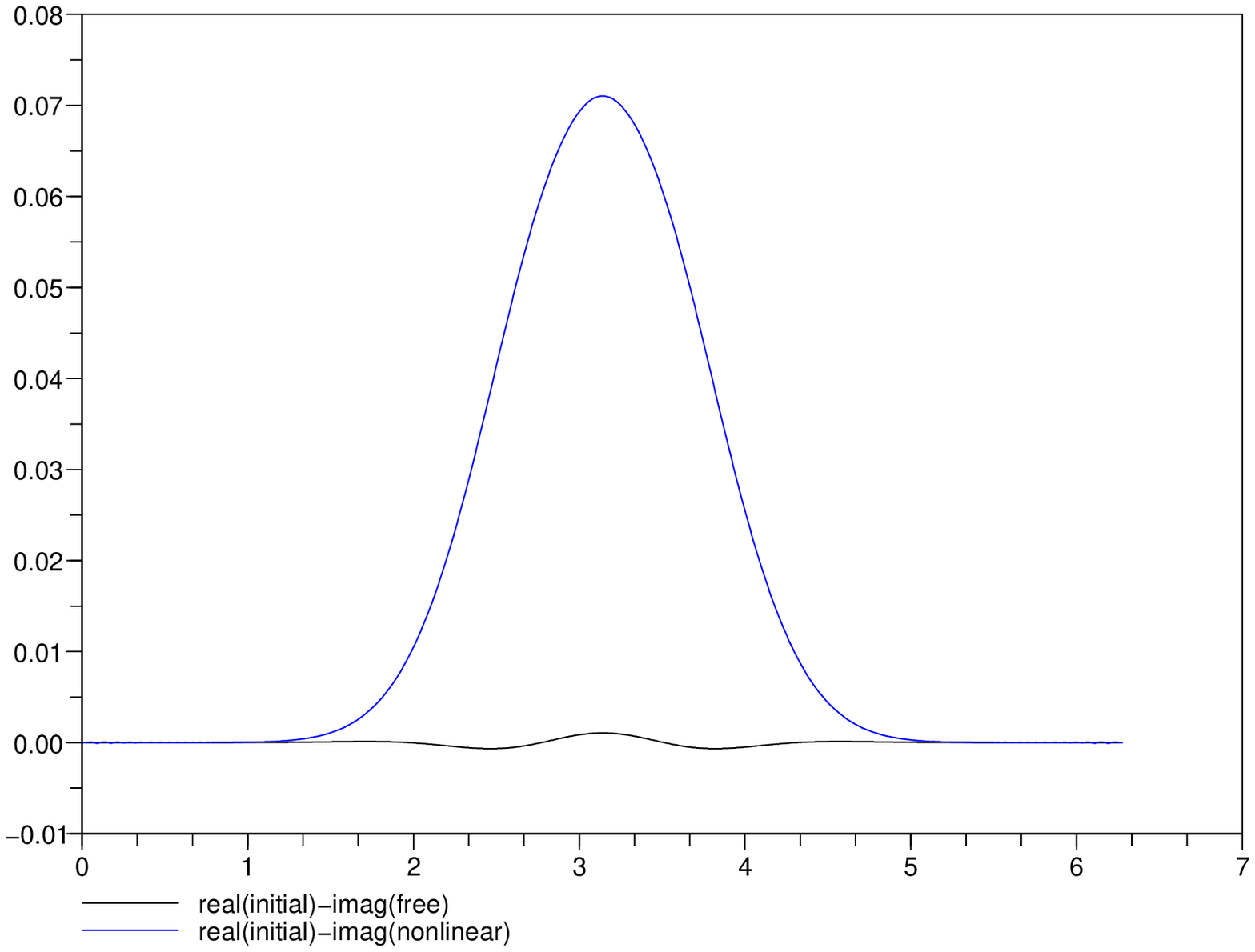,width=0.45\linewidth}
\epsfig{file=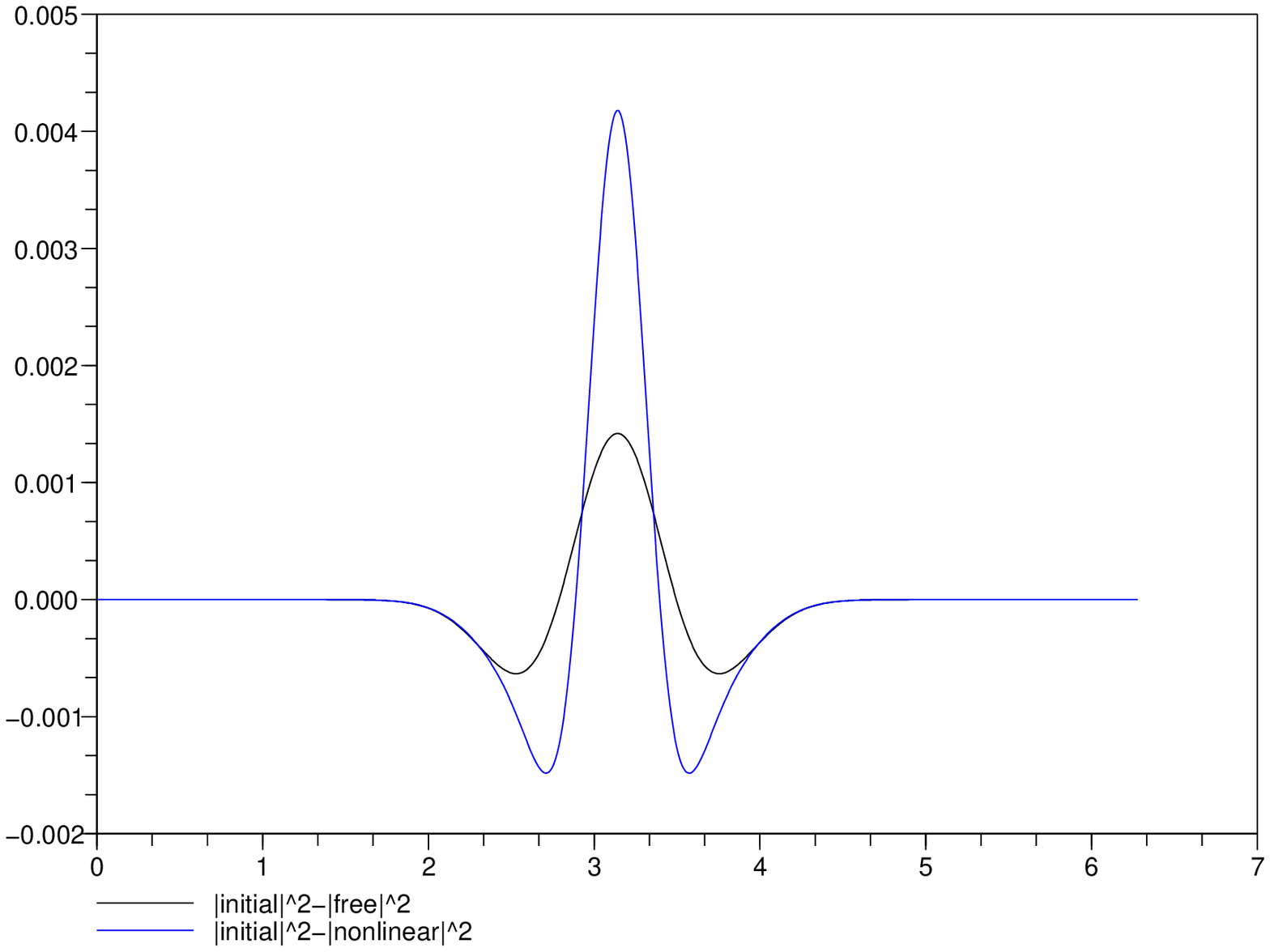,width=0.45\linewidth}}
\caption{Same as Fig.~\ref{pf1}, but $\sigma=2$ and $\alpha=2$. (Critical
case)} \label{pf2} 
\end{figure}

We now put $\sigma=\alpha=2$ and the outcome is displayed on
Fig.~\ref{pf2}; we still compare the same quantities. Absolute errors
on wave functions (left side) are much bigger for $u^\e$ in this
case. 
In particular, no new frequencies appear in the numerical
solutions. The nonlinear effect boils down to a small change on the
modulus of $u^\e$. 

\subsection{Supercritical case}

\begin{figure}[ht]
\centerline{\epsfig{file=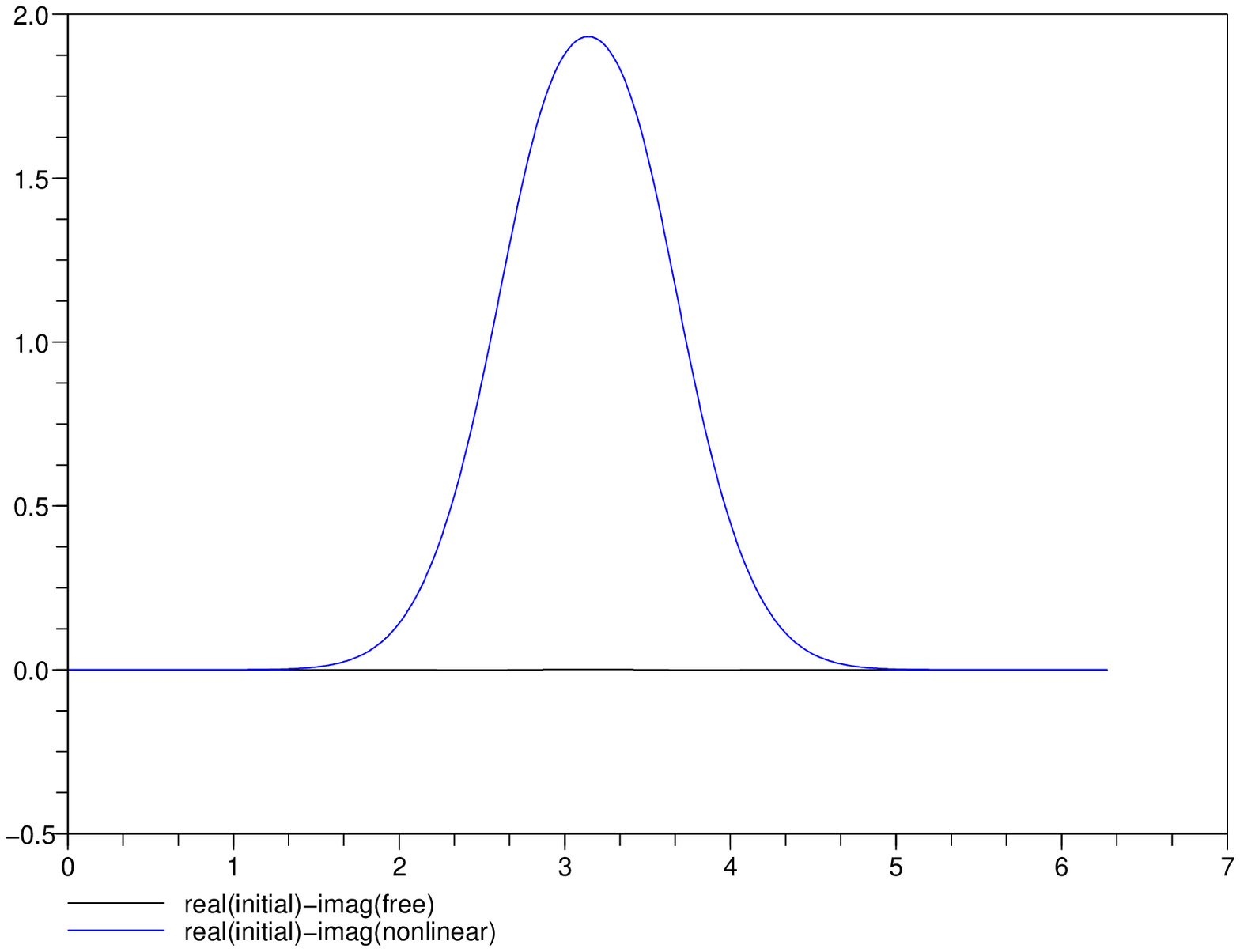,width=0.45\linewidth}
\epsfig{file=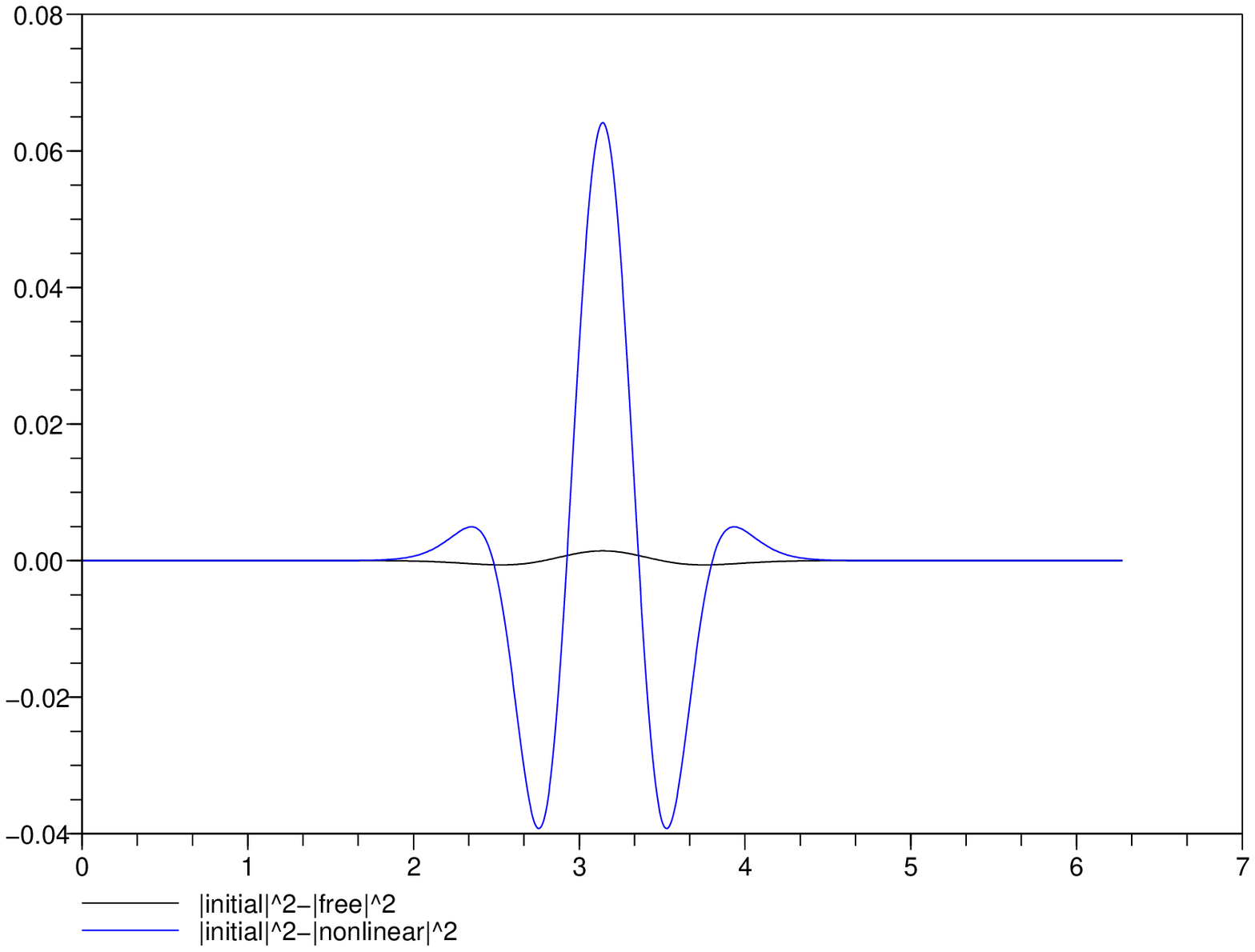,width=0.45\linewidth}}
\caption{Same as Fig.~\ref{pf1}, but $\sigma=2$ and
$\alpha=1.5$. (Supercritical case)} \label{pf3} 
\end{figure}

We close this first series of tests by considering  $\sigma=2$ and
$\alpha=1.5$ as shown in Fig.~\ref{pf3}. Of course, as no pointwise
convergence is expected in this case, absolute errors are even bigger
for both wave functions (left side) and moduli (right
side).  
%
Of course, the size of 
the error on the modulus is much bigger too, and one should be
extremely careful about the credit to give to the numerical
simulations in the supercritical case. Indeed, this is a regime where
a small error can be amplified at leading order (see \cite{CaCascade,CaARMA}). 

\section{Visualization of the scattering operator}\label{sec:scattnum}
\vspace*{-0.5pt}
\noindent

We aim now at illustrating the results on scattering theory through
numerical computations still achieved through time-splitting FFT
schemes. The algorithm we used for the approximation of
the scattering operator is based on a nonlinear time-splitting routine
flanked by two free evolution steps (implemented the way recalled in
the previous section): 
$$
S \psi \simeq U(-T)\circ U_{NL}(2T) \circ U(-T) \psi, \qquad T \gg 1,
$$
with $U$, $U_{NL}$ standing for the solution operators of equations
\eqref{eq:schrodlibre} and \eqref{eq:nls} in 1-D with $\e =1$
respectively. We used $T=55$ in the computations hereafter. 
 
As no small parameter $\e$ is present in the problem, one may think
that no major obstacle exists in carrying out this program; this isn't
correct as the free evolutions can (and do!) dramatically increase the
size of the computational domain for large $T$. It is interesting to
notice that, in case one wants to use FFT-based schemes, both the
direct computation for small $\e$ and the scattering operator
approximation lead to a ``large computational domain difficulty": in
the Fourier space for the first case, in the usual space for the
second. 

A way to understand the scattering operator is to
visualize the average effects of the nonlinearities appearing in
equations of the form 
\begin{equation}
  \label{eq:nls2}
  i \d_t u_\lambda +\frac{1}{2}\partial_x^2 u_\lambda
  =\lambda |u_\lambda|^{2\si}u_\lambda, \qquad u_\lambda(t=0,x)=\exp(-5 x^2),
\end{equation}
for various values of $\sigma \geq 1$ and $\lambda$. Intuitively, as
 $\sigma$ increases, the nonlinearity becomes shorter
 range. Similarly, as $\lambda$ increases, the nonlinearity becomes
 stronger, and it should take a larger amount of time before we can
 consider it has become negligible. In all the tests we performed, it
 was somehow surprising to observe how fast the algorithm converges:
 one does not have to consider ``very large'' values of $T$ so that  
\begin{equation*}
U(-T)\circ U_{NL}(2T) \circ U(-T) \psi
\end{equation*}
becomes stable and visually independent of $T$. 

\subsection{Quintic nonlinearity ($\sigma=2$)}

\begin{figure}[ht]
\centerline{\epsfig{file=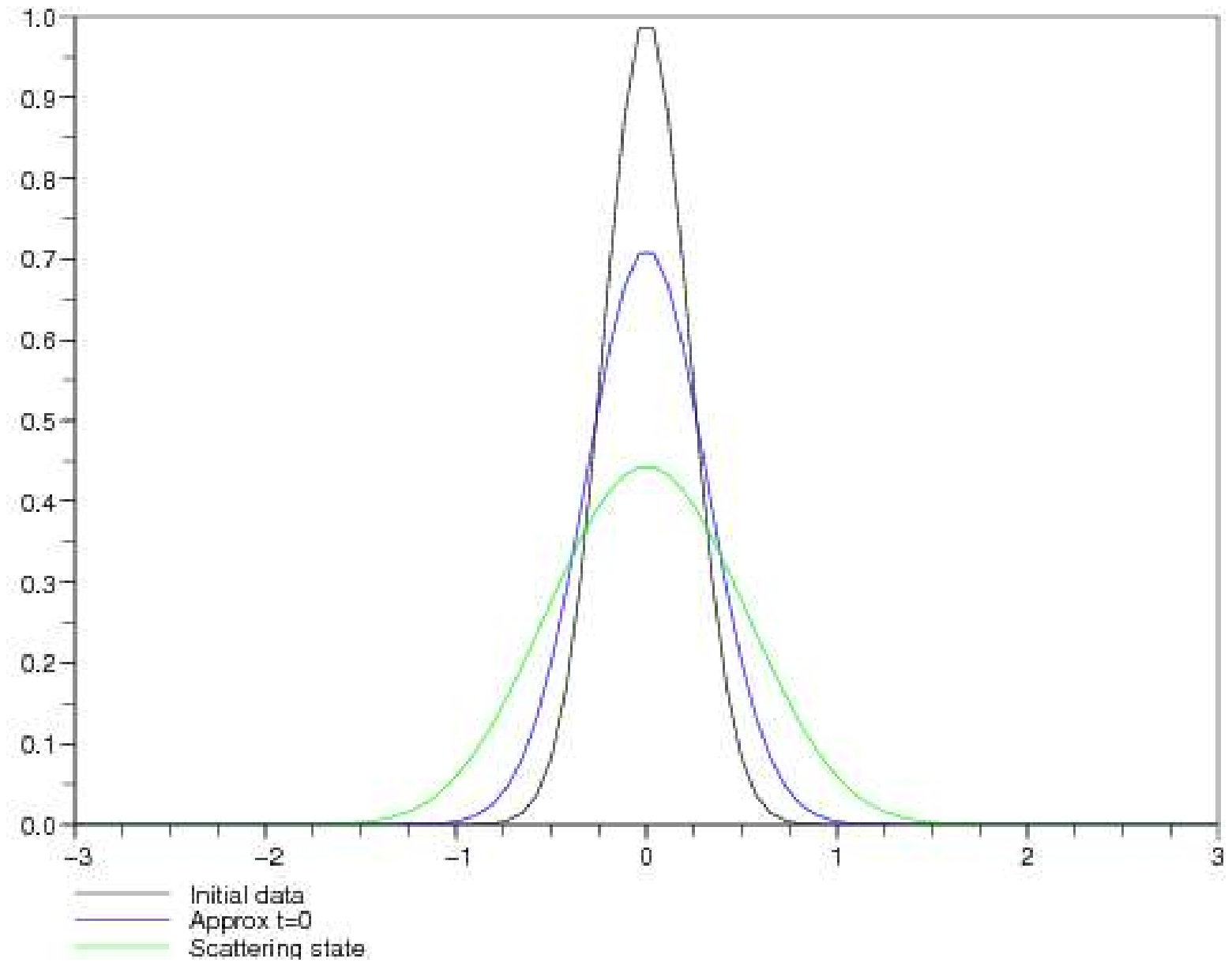,width=0.45\linewidth}  
\epsfig{file=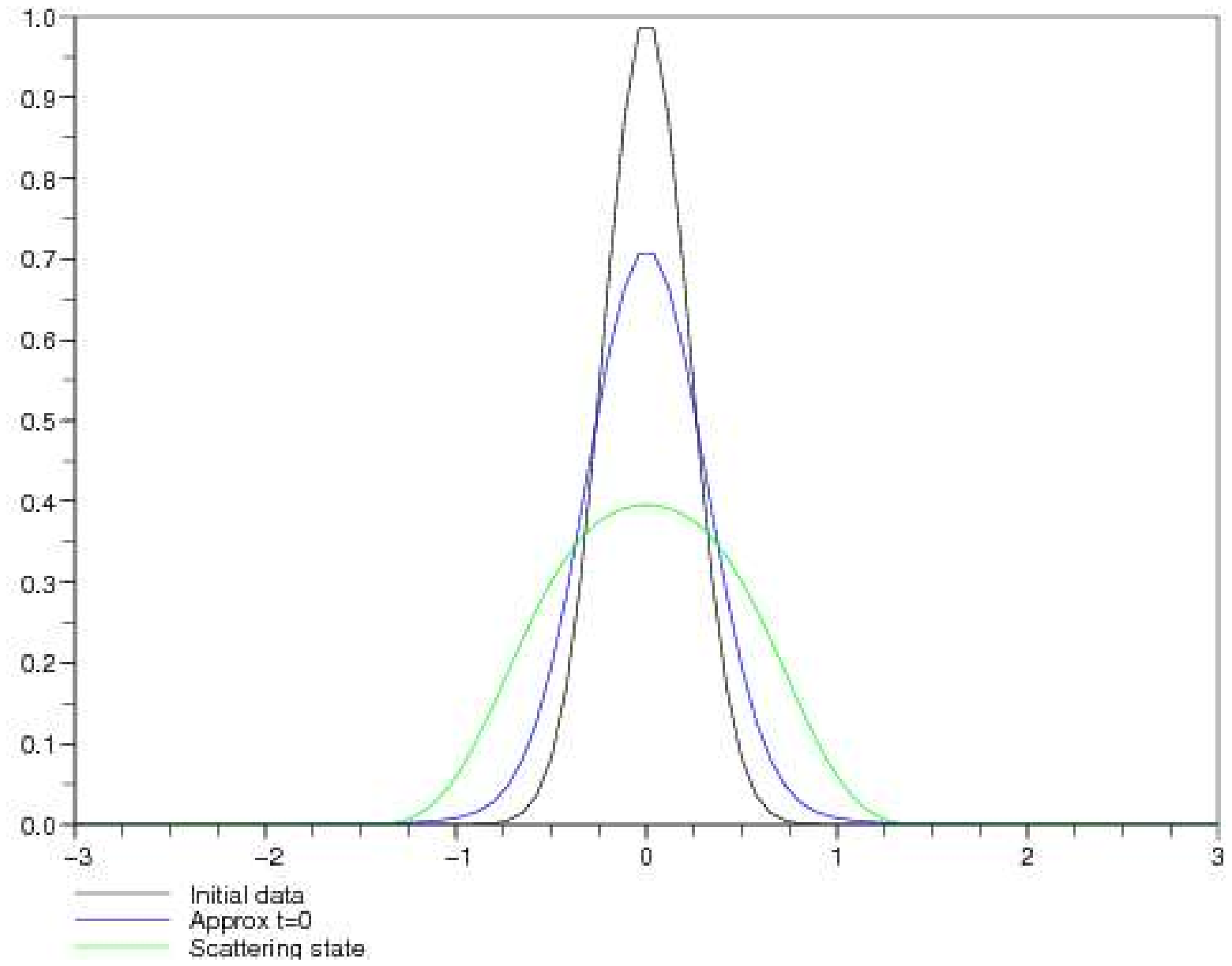,width=0.45\linewidth}}
\caption{Initial data, numerical solution at $T=0$, scattering state
  for \eqref{eq:nls2} with 
$\sigma=2$ and $\lambda=1$ (left), $\lambda=5$ (right)} \label{run1}
\end{figure}

\begin{figure}[ht]
\centerline{\epsfig{file=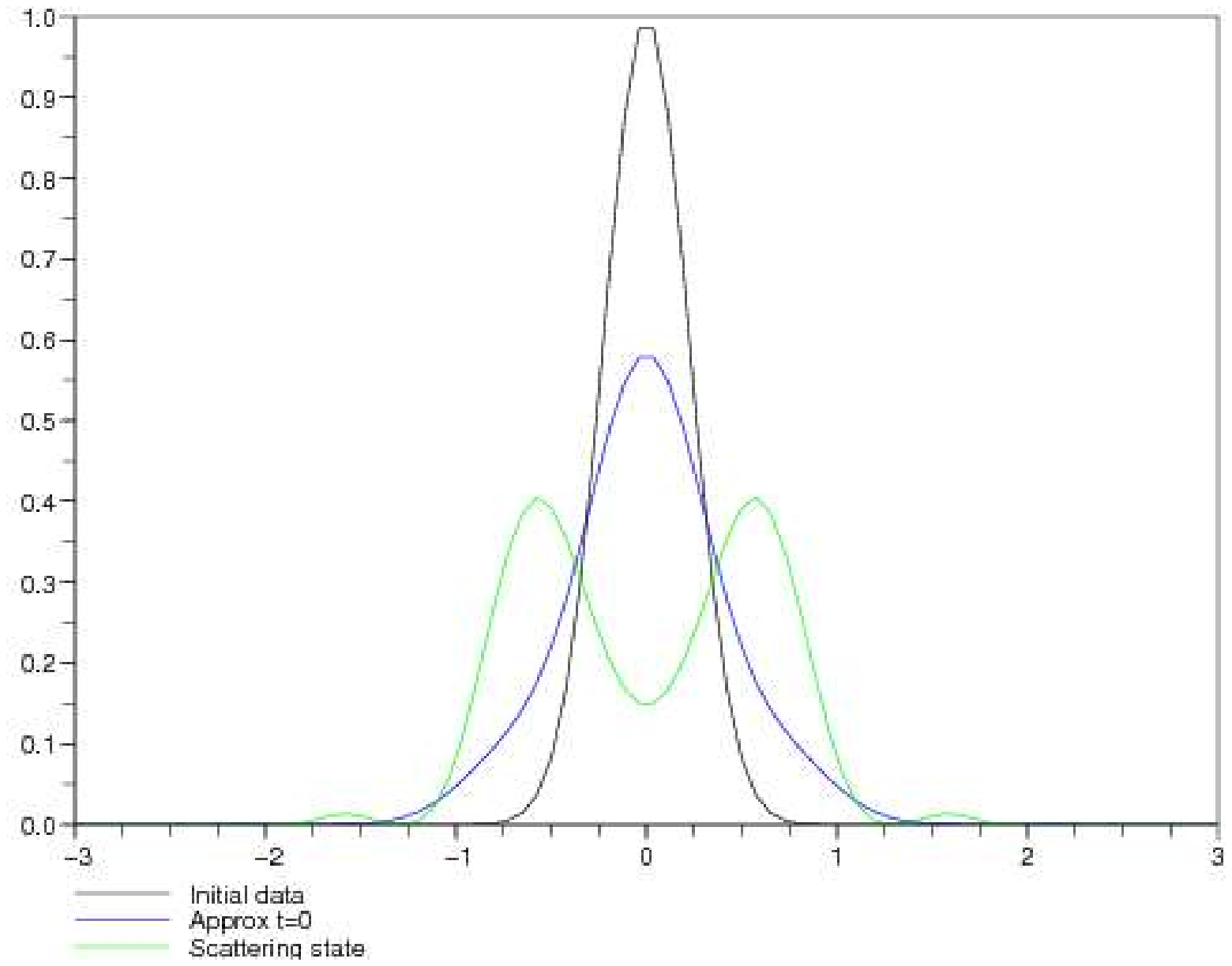,width=0.45\linewidth} 
\epsfig{file=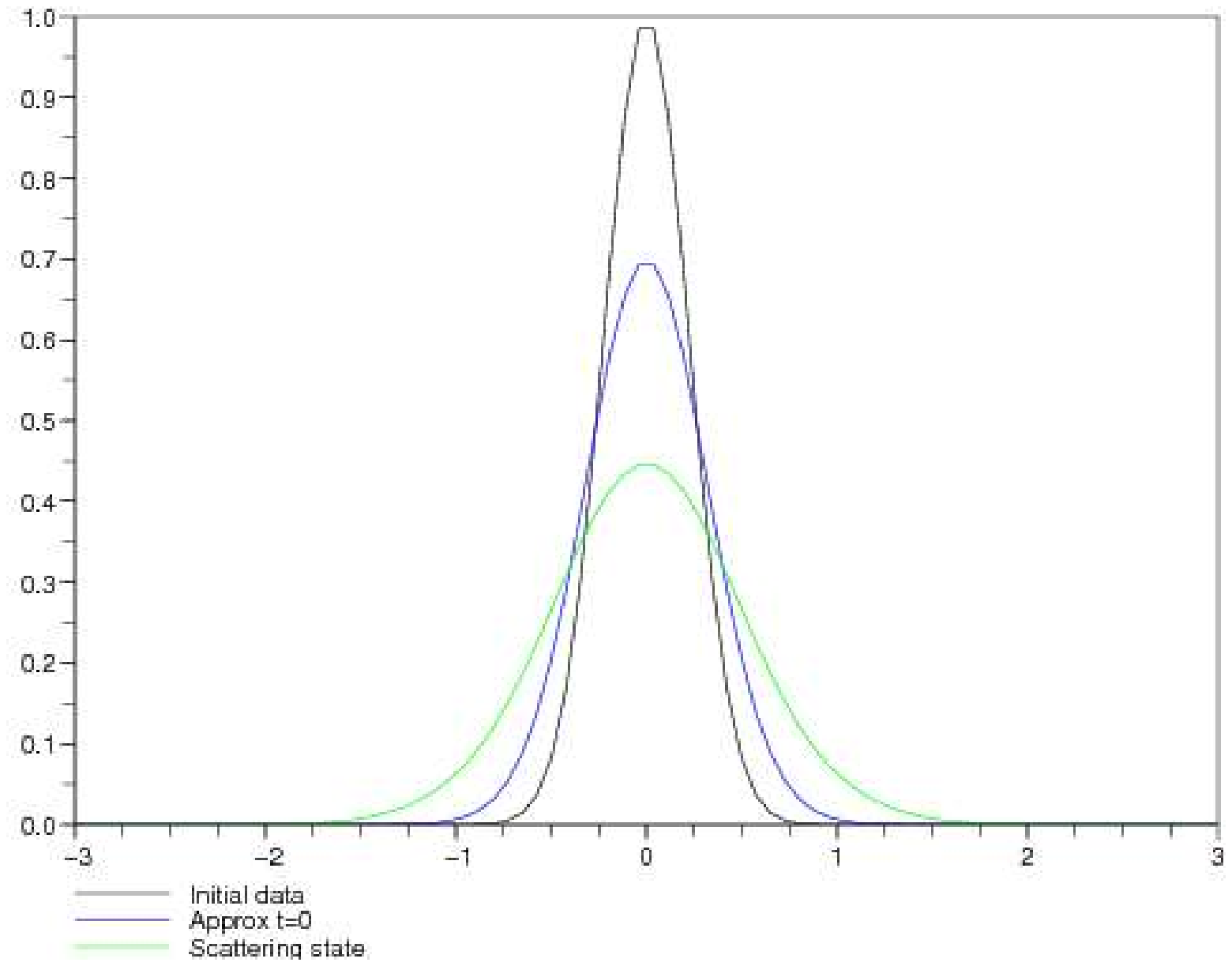,width=0.45\linewidth}}
\caption{Same as Fig.~\ref{run1}, but $\lambda=25$ (left),
  $\lambda=-1$ (right)} \label{run2} 
\end{figure}

The parameter $\lambda$ controls in some sense the strength of the
nonlinearity\footnote{through the stiffness of the associated
  differential equation.} inside \eqref{eq:nls2}, as can be seen on
Fig.~\ref{run1}. This figure displays the position density of the
initial data, the scattered solution for $T=55$ and a ``mixed state"
$\tilde u_\lambda(t=0)=U(-T)U_{NL}(T)u_\lambda(t=0)$. As our
time-splitting/FFT algorithm preserves only the $L^2$ norm, but not
the Hamiltonian, we first restricted ourselves to moderate values of
$\lambda \geq 0$ (defocusing case).  
However, as a numerical experiment, we wanted to display the outcome
of our scheme for the stronger case $\lambda=25$ on Fig.~\ref{run2}:
notice the change of shape in the scattered solution. Moreover, on
this figure, we also tried to show what happens for $\lambda=-1$, that
is to say for the focusing case despite there may be finite time
blow-up (but there is scattering for small data). We checked that 
the energy associated to this data is (and remains) 
positive, a case where the virial identity , \cite{CazCourant},
does not imply blow-up. The computational domain for these runs 
was $[-100\pi, 100\pi]$ and $2^{13}-1$ Fourier modes were used. 
\subsection{Power 3 ($\sigma=1.5$)}

\begin{figure}[ht]
\centerline{\epsfig{file=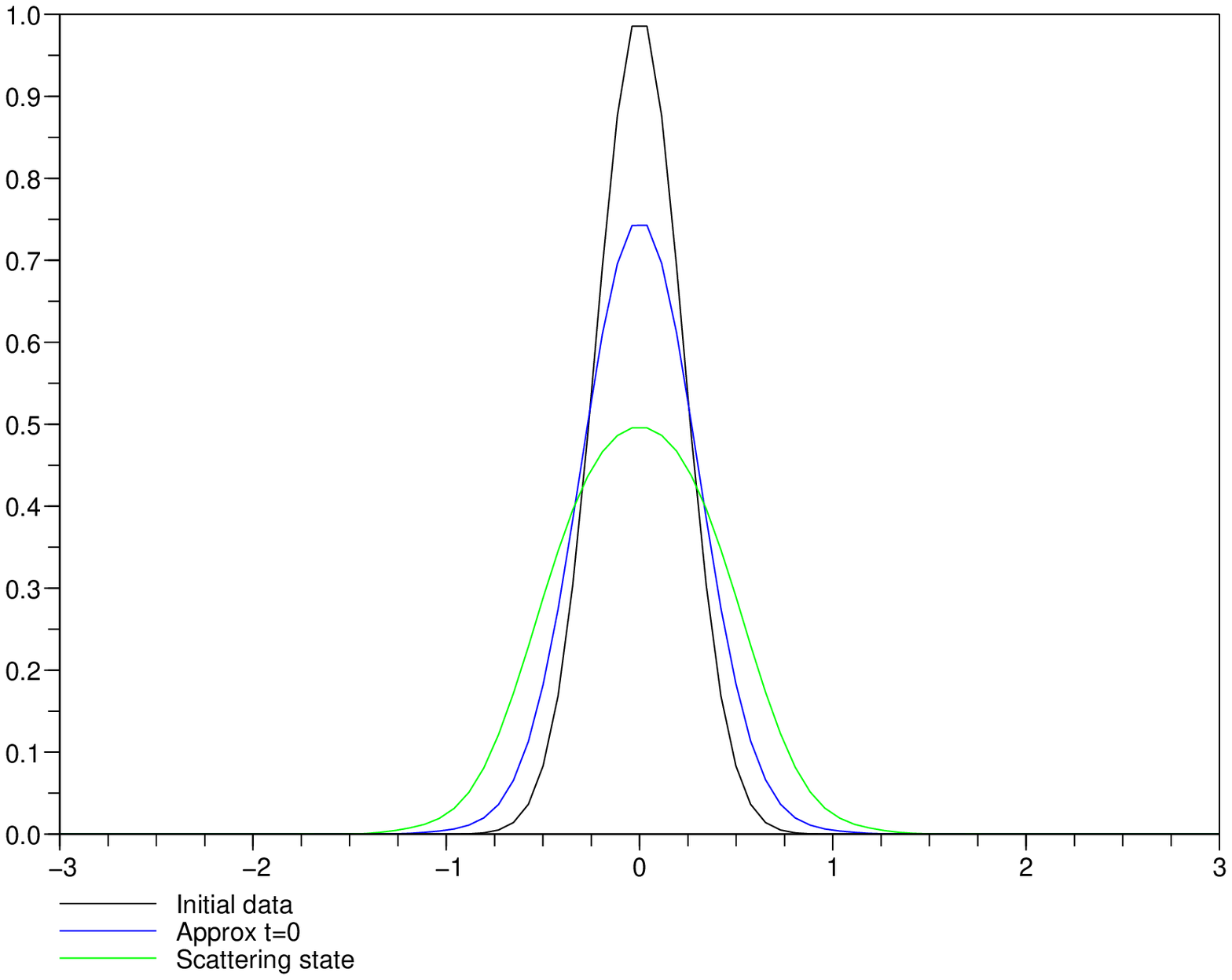,width=0.45\linewidth}  
\epsfig{file=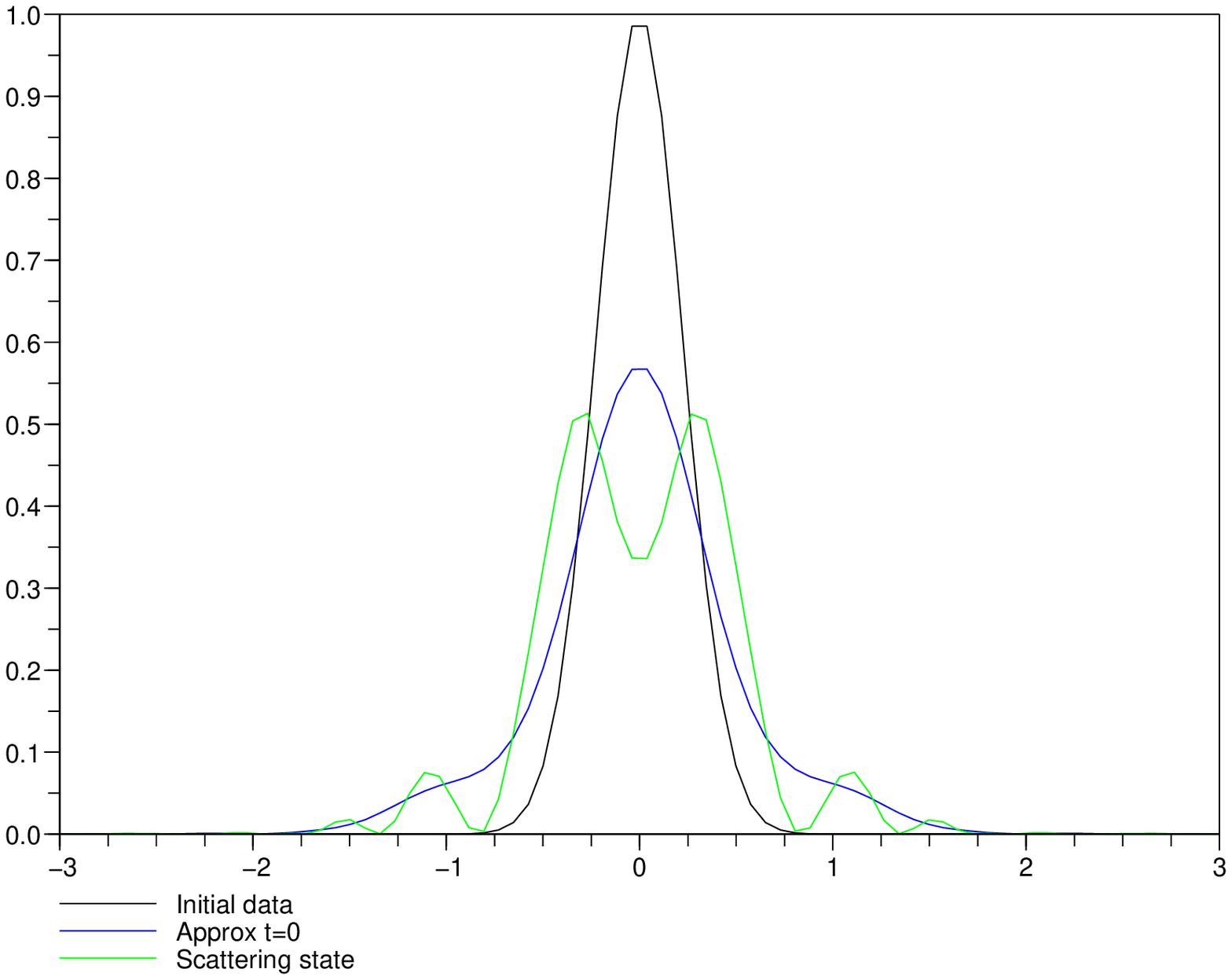,width=0.45\linewidth}}
\caption{Same as Fig.~\ref{run1}, but $\sigma=1.5$ and $\lambda=5$
  (left), $\lambda=25$ (right)} \label{run3} 
\end{figure}

Now let's observe the effects of lowering the $\sigma$ value while
keeping other parameters equal, see Fig.~\ref{run3}. It is interesting
to see that the change of shape appearing for $\lambda=25$ is stronger
than in the preceding case. On the contrary, the increase of the
numerical solution's support is slightly less important. 
This hints that increasing the $\sigma$ value tends to expand the
support of the scattered solution whereas increasing the $\lambda \geq
0$ value (defocusing case) leads to an oscillatory behavior. However,
we stress that since the energy, 
$$
E(t):=\frac 1 2 \|\d_x u^\lambda (t) \|_{L^2}^2 + 
\frac \lambda {\sigma + 1}\| u^\lambda(t)\|_{L^{2\si +2}}^{2\si +2}=E(0),
$$
of the numerical solution
changes more with a bigger $\lambda$ (its mass being always kept
constant), these oscillations might be spurious. We actually don't
know how this fact can be decided; our profiles have been checked to
be stable on a finer grid. 

\subsection{Power 6 ($\sigma=3$)}

\begin{figure}[ht]
\centerline{\epsfig{file=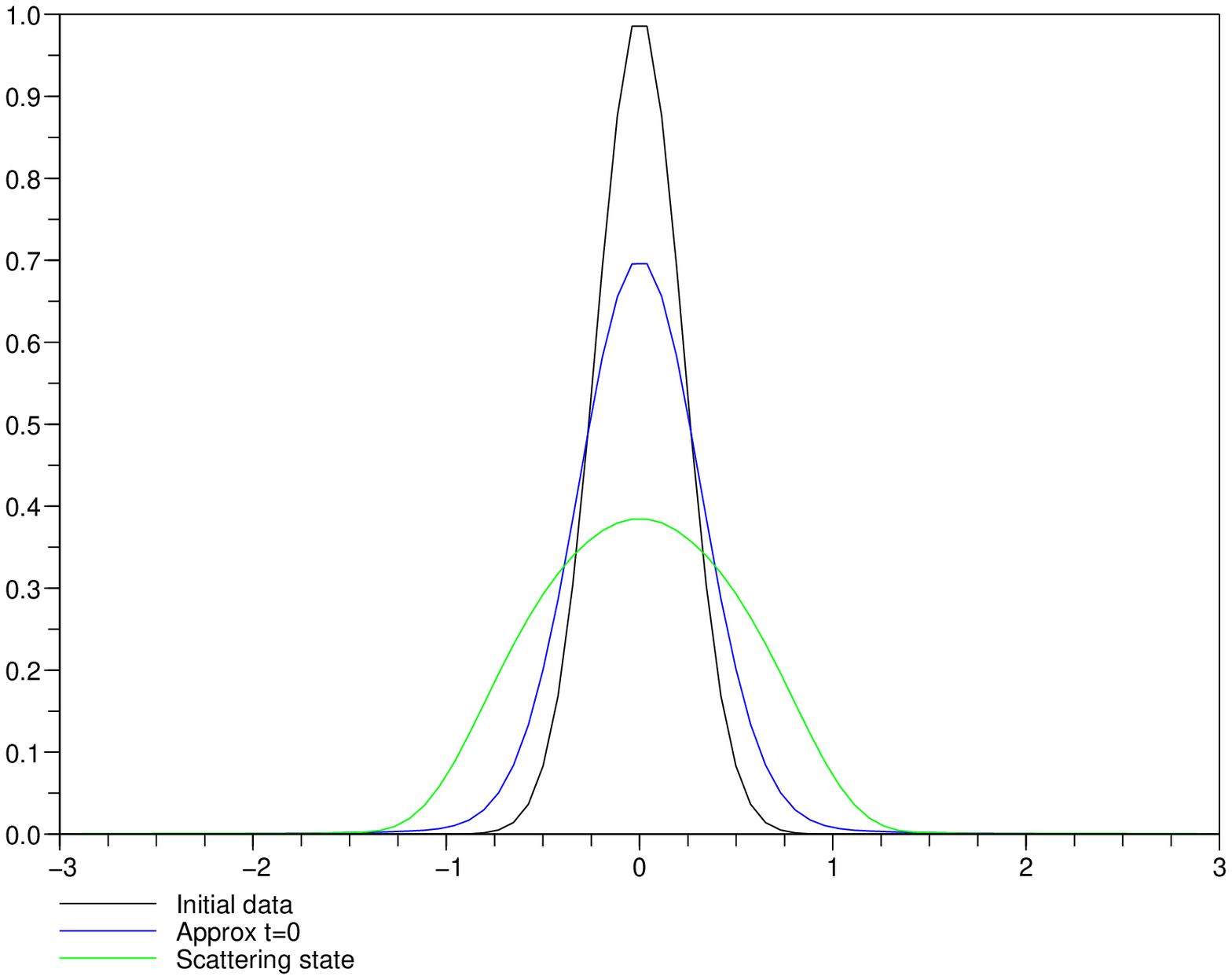,width=0.45\linewidth}
\epsfig{file=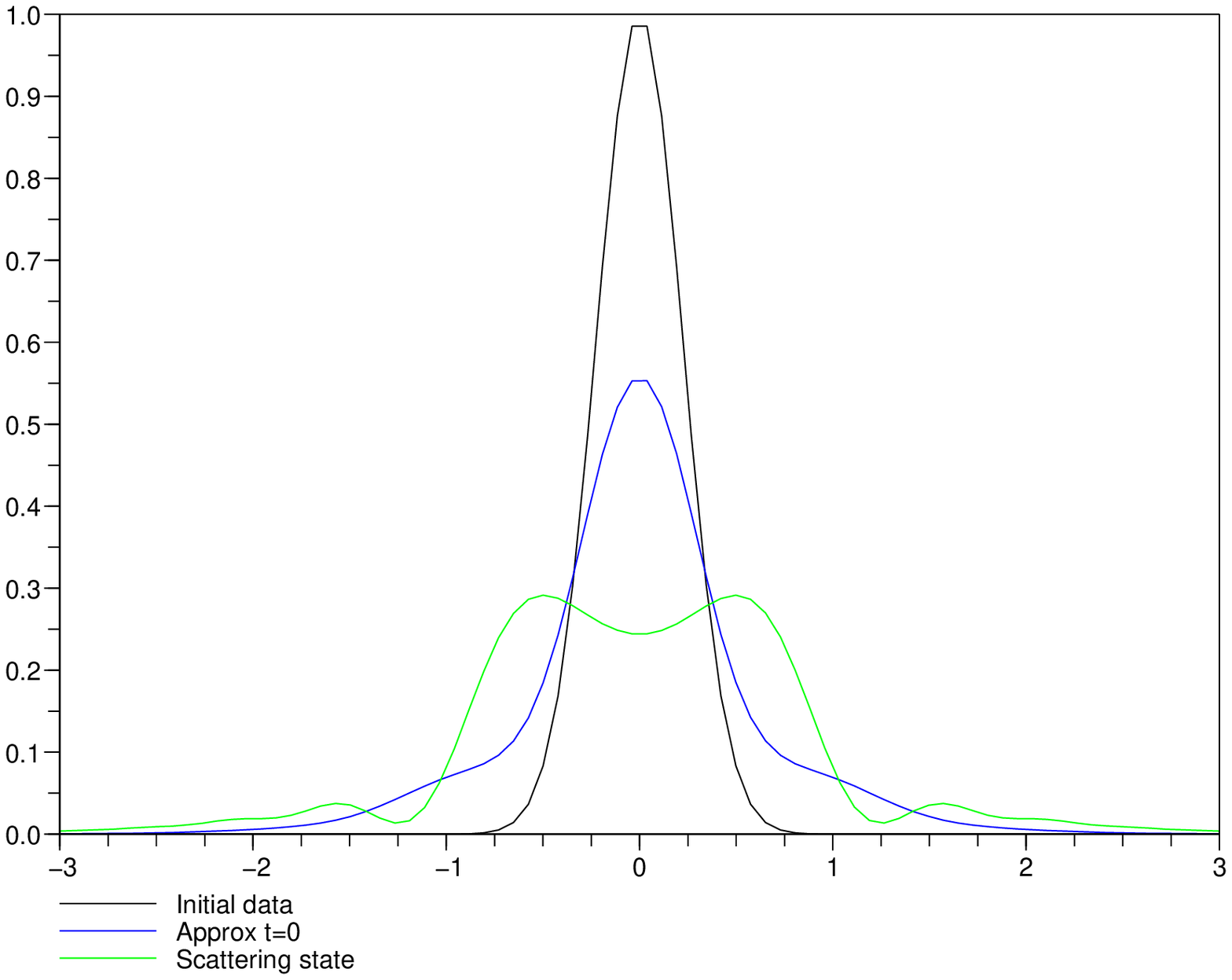,width=0.45\linewidth}}
\caption{Same as Fig.~\ref{run1}, but $\sigma=3$ and $\lambda=5$
  (left), $\lambda=25$ (right)} \label{run4} 
\end{figure}

In order to get some numerical evidence about the dependence of the
scattered solution on $\sigma$, we display on Fig.~\ref{run4} the
outcome for $\sigma=3$. It is quite clear that the scattered solutions
for both values of $\lambda$ are less peaked. Their support is bigger
and the oscillations for $\lambda=25$ are weaker, their frequency
remained the same though. 
This agrees with the behavior we sketched in the preceding subsection
as $\sigma$ and $\lambda$ vary. 

\section{Experiments on a cusp caustic}\label{sec:cusp}
\vspace*{-0.5pt}
\noindent

\begin{figure}[p]
\centerline{\epsfig{file=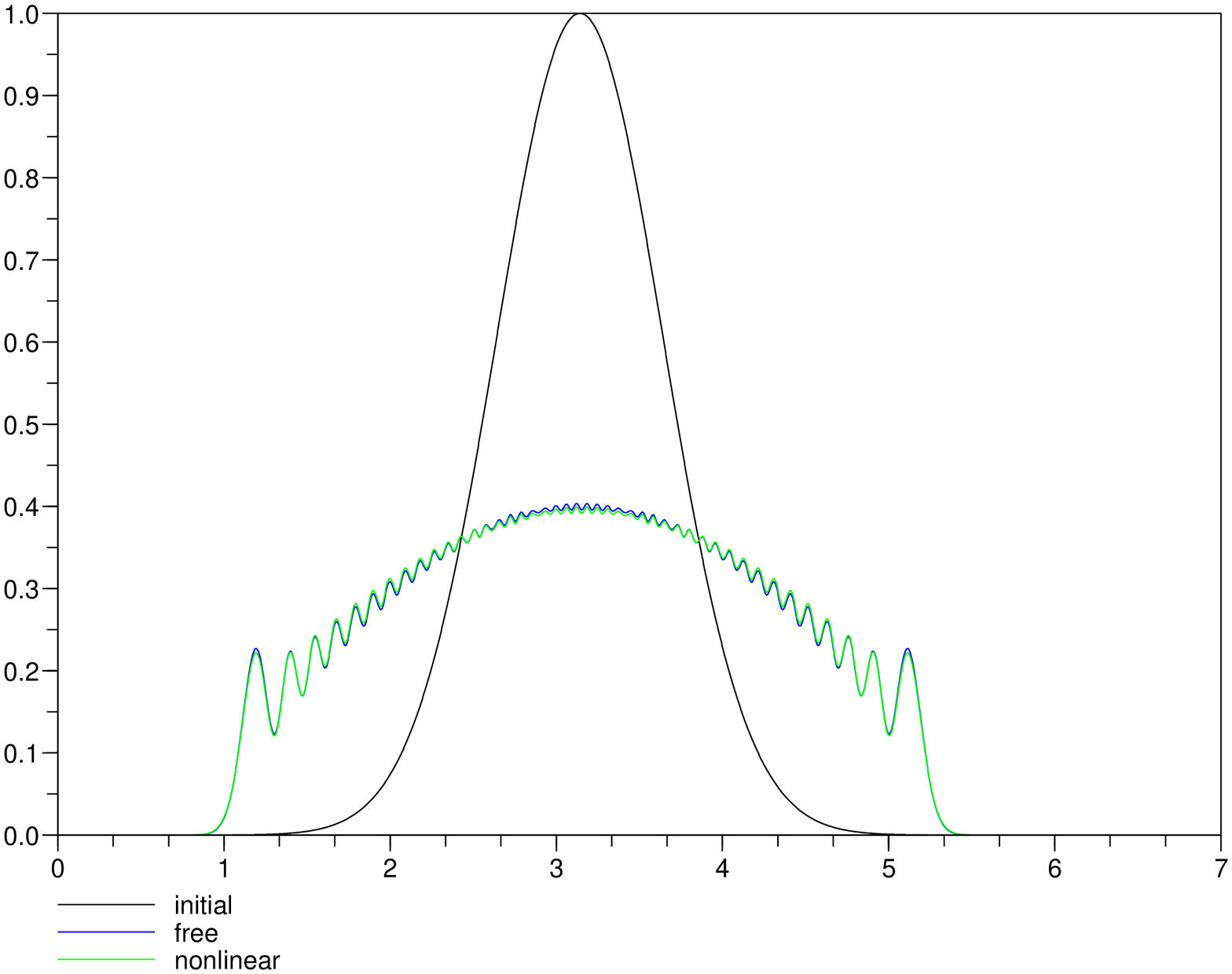,angle=270,width=0.7\linewidth}} 
\centerline{\epsfig{file=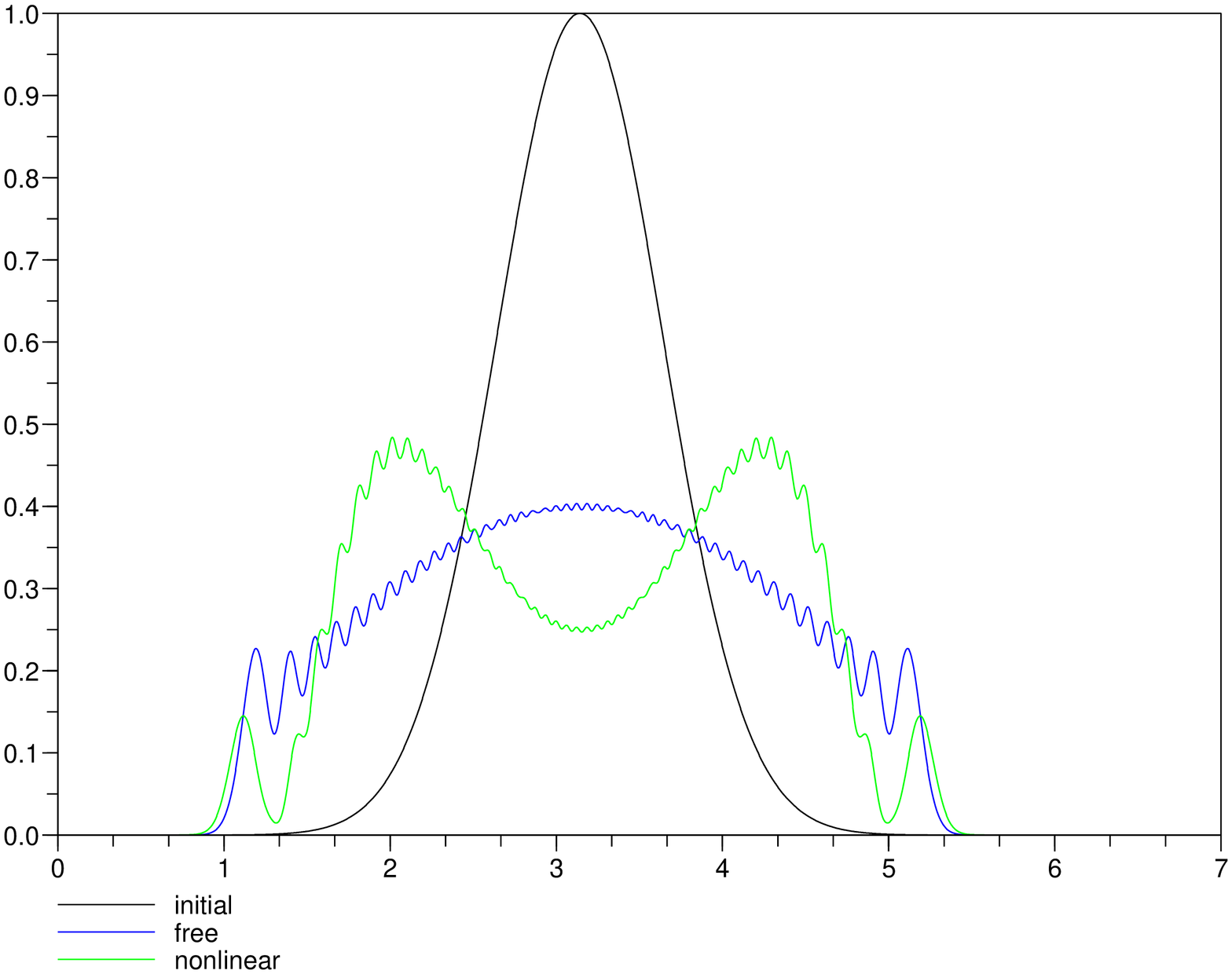,angle=270,width=0.7\linewidth}}
\centerline{\epsfig{file=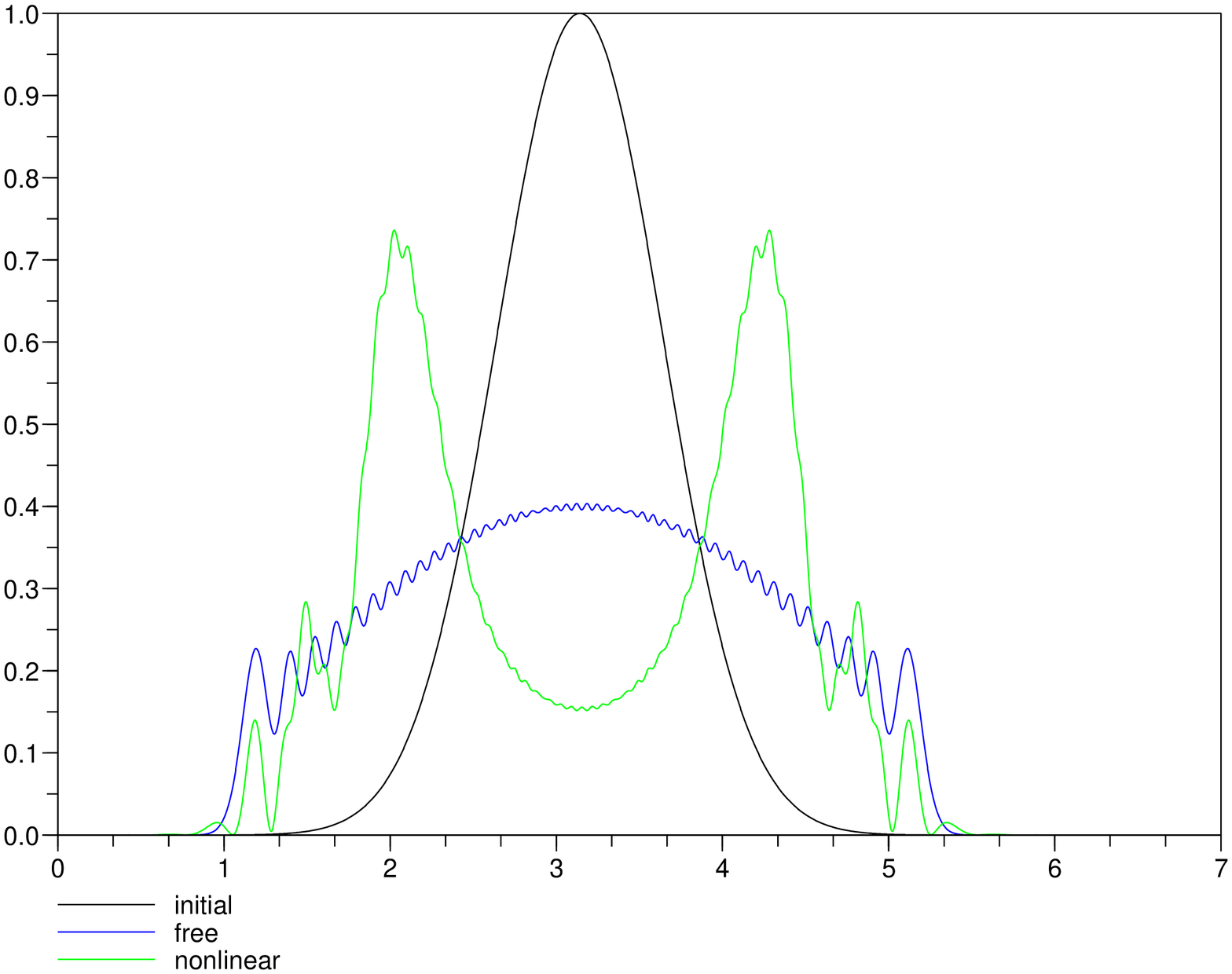,angle=270,width=0.7\linewidth}}
\caption{Position densities in the cusp caustic: $\alpha=4,3,2$ (top to bottom)} \label{run-cusp}
\end{figure}

Let us now go back to comparing the quadratic observables generated by
numerical approximation of equations \eqref{eq:nls} and
\eqref{eq:schrodlibre} endowed with a small parameter $\e$ in 1-D. In
this section we fixed $\e = 1/150$. Figure \ref{run-cusp} displays the
position density of the initial data for both equations, i.e. 
$$
u^\e(t=0,x)=v^\e(t=0,x)=\exp\Big(-2(x-\pi)^2-i\cos(x)/\e \Big), \qquad
x \in [0,2\pi], 
$$
together with the position density of the numerical approximations of
\eqref{eq:nls}, \eqref{eq:schrodlibre} in $T=3.5$. The point here is
to investigate what happens for the case of such a self-interfering
Gaussian pulse, since no scattering theory is known for this
problem. What we would like to check is whether the theoretical
results on the focus point recalled and visualized in the preceding
sections can be thought of as a guideline for this more complex case
involving a non-trivial caustic. 
We shall observe position densities for the unique value of $\sigma=4$
as a similar behavior has been seen to hold for different
nonlinearities with convenient values of $\alpha$. $4095$ Fourier
modes have been used in order to produce these results. 

\subsection{Subcritical picture: $\alpha=4$}

This case could be referred to as subcritical since it is noticeable on
the top of Fig.~\ref{run-cusp} that the free and the nonlinear
numerical solutions do agree for this reasonably small value of
$\e$. In particular, the frequencies of oscillations are
identical. This is very similar compared to the behavior investigated
in \cite{CaIUMJ}. 

\subsection{Critical picture: $\alpha=3$}

\begin{figure}[ht]
\centerline{\epsfig{file=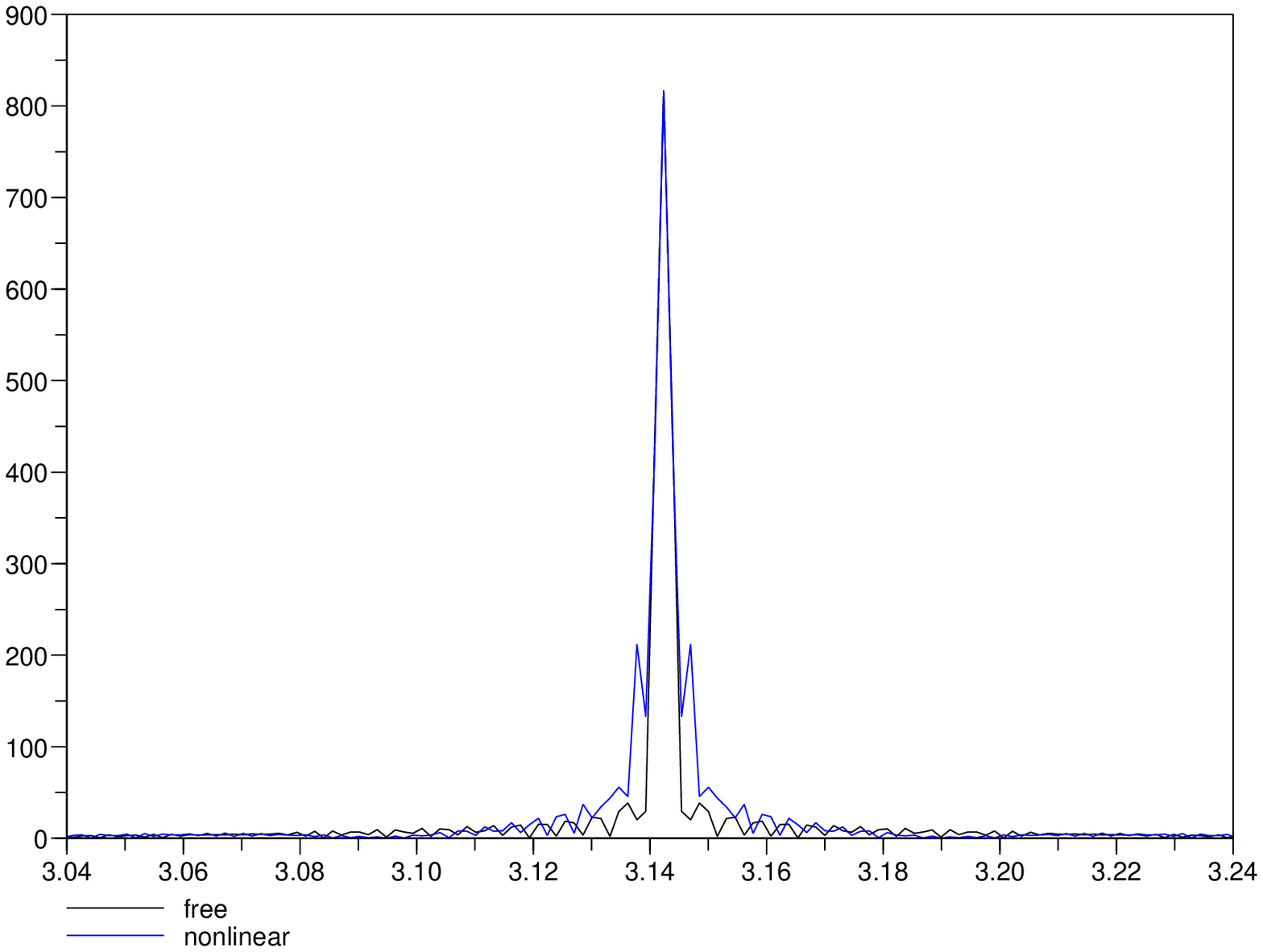,width=0.45\linewidth}
\epsfig{file=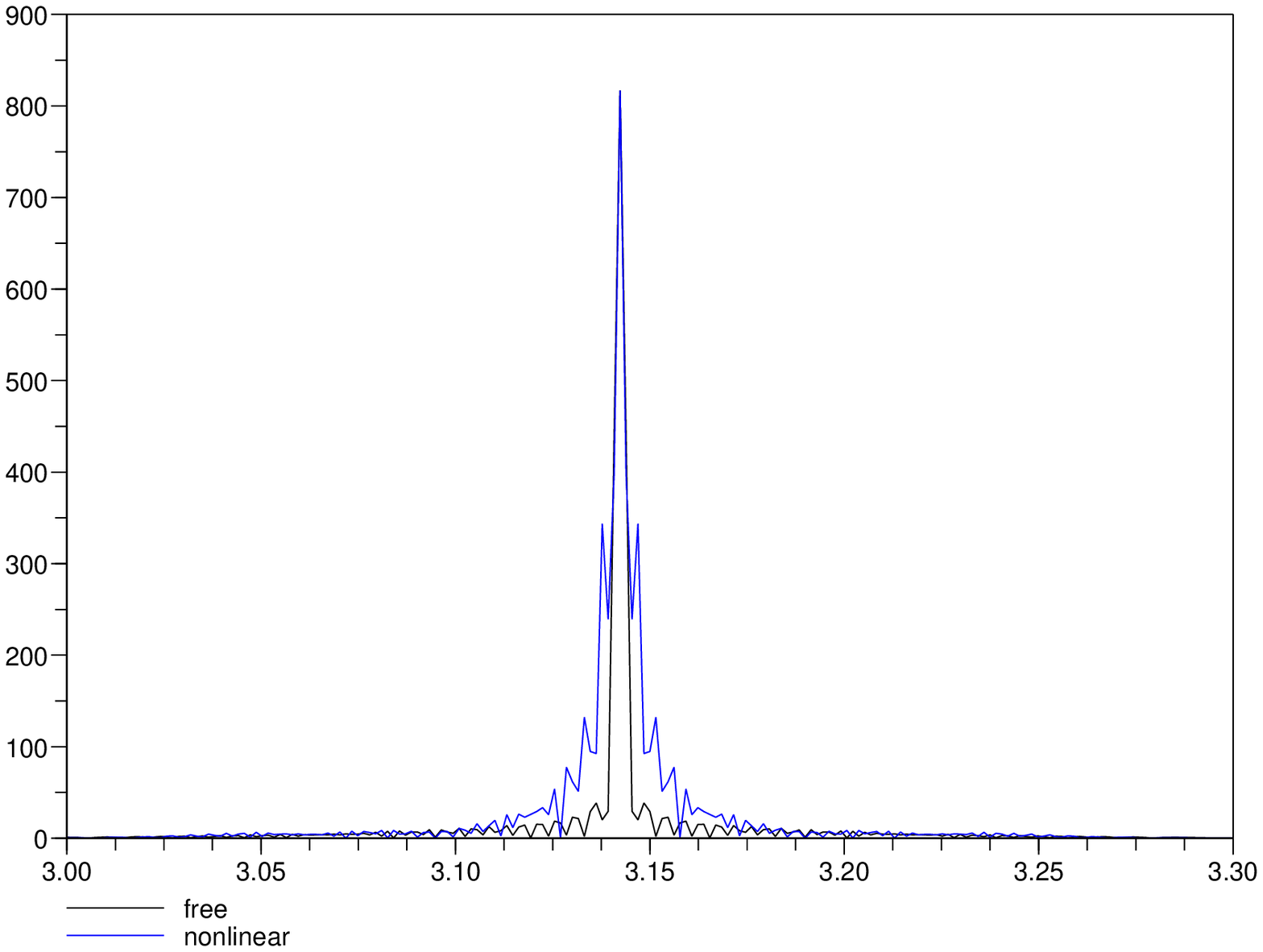,width=0.45\linewidth}}
\caption{Fourier transforms of the position densities with $\alpha=3$
  (left), $\alpha=2$ (right)} \label{fft-cusp} 
\end{figure}

The parameter $\alpha$ is now in a ``critical range" as we observe
that both solutions differ much more, but the frequency of the
oscillations looks like being still the same in both cases. In order
to establish this fact, we display on the left of Fig.\ref{fft-cusp}
the FFT of the position densities: a peak at the same frequency is
clearly noticeable.  
The nonlinear effect manifests itself through a change of order zero
in the moduli, as we already observed on the right side of
Fig.~\ref{pf2}; notice also the similarity with the scattering state
shown on Fig.~\ref{run2} (right). This does agree with the $\alpha_c$ 
value derived in Section \ref{sec:heur} 

\subsection{Supercritical picture: $\alpha=2$}

In this last case, there is no similarity no more between the
approximate solutions of \eqref{eq:nls}, \eqref{eq:schrodlibre}, as
seen on both Fig.~\ref{run-cusp} and \ref{fft-cusp}. Especially, the
right side of Fig.~\ref{fft-cusp} reveals that new frequencies show up
inside the position density of the nonlinear solution. We have
therefore a change of order zero in the moduli and in the
frequency. This is of course reminiscent of Fig.~\ref{pf3} in which a
frequency doubling seems to show up in the supercritical regime. 

\section{Conclusion}\label{sec:concl}
\vspace*{-0.5pt}
\noindent

We have presented the semi-classical limit for the nonlinear
Schr\"odinger equation in the presence of a caustic. When the caustic
is reduced to a point, the numerical experiments are in good agreement
with the analytical results as far as the notion of criticality is
concerned. However in the critical case, described by a nonlinear
scattering operator, the leading order nonlinear effects are rather
hard to visualize in the semi-classical limit. This is why we
simulated the scattering operator in a separate way. 

Our numerical
tests give encouraging evidence of new phenomena concerning the
phase of the wave in the supercritical case when a focal point is
formed (appearance of new frequencies). In the presence of a cusp caustic, 
the numerical experiments are in
good agreement with the heuristic arguments that we presented here, for
which no rigorous justification is available so far.

\end{document}